\documentclass[preprint,12pt]{elsarticle}

\usepackage{amssymb}
\usepackage{amsmath}
 \usepackage{amsmath,amssymb,amsfonts,amsthm}
 \usepackage{graphicx}
\usepackage{layout}
\usepackage[latin1]{inputenc}
\usepackage[french,english]{babel}
\usepackage[french]{babel}
\usepackage{venndiagram}
\usepackage{cancel}
\newtheorem{theorem}{Theorem}
\newtheorem{definition}{Definition}

\newtheorem{proposition}{Proposition}
\newtheorem{lemma}{Lemma}

\newtheorem{example}{Example}
\newtheorem{remark}{Remark} 
\usepackage{enumitem}

\newcommand{\vertiii}[1]{{\left\vert\kern-0.25ex\left\vert\kern-0.25ex\left\vert #1 
    \right\vert\kern-0.25ex\right\vert\kern-0.25ex\right\vert}} 

\journal{ArXiv}

\begin{document}

\begin{frontmatter}



\title{A New Crossnorm That Preserves Unconditional Bases in Banach Spaces}
 \author[label1]{Rafik Karkri}
 \affiliation[label1]{organization={FSK, Ibn Tofail University},
             addressline={P.O.Box 133},
             city={Kenitra},
             postcode={14000},
             country={Morocco}}
\ead{karkri.rafik@gmail.com}
\author[label1]{Samir Kabbaj}
 \ead{samkabbaj@yahoo.fr}
\begin{abstract}
Let $\alpha$ be a tensor norm (i.e., a uniform reasonable crossnorm) on the class of all algebraic tensor products of Banach spaces $E \otimes F$. We say that $\alpha$ preserves unconditionality if, for every pair of Banach spaces $E$ and $F$ with  unconditional Schauder bases (USBs), the completion $E \otimes_{\alpha} F$ also admits a USB.

It is well known that none of Grothendieck's fourteen natural tensor norms satisfy this unconditionality-preserving condition. Moreover, the existence of a tensor norm $\alpha$ with this property remains an open question.

In this paper, we construct for every such pair $(E,F)$ a new reasonable crossnorm $\alpha$. This norm has the surprising property that---despite being generally non-uniform---the space $E \otimes_{\alpha} F$ nevertheless admits a USB.
\end{abstract}

\end{frontmatter}



\section{Introduction}
In 1961, Gelbaum and Gil de Lamadrid \cite[Theorem 5, p. 1285]{Gelbaum.Lamadrid.1961} observed that the canonical Schauder basis of the space $\ell_{2}\otimes_{\varepsilon}\ell_{2}$ is not unconditional.
In 1970, Kwapie\'{n} and Pe\l czy\'{n}ski \cite[Corollary 3.1, p. 58]{Kwapien.Pelski.1970} proved  that
neither $\ell_{p}\otimes_{\pi}\ell_{q^{*}}$ nor $\ell_{p^{*}}\otimes_{\varepsilon}\ell_{q}$ (for $p\geq q>1$) admits an unconditional Schauder basis. These results demonstrate that both the injective $(\varepsilon )$ and projective $(\pi)$ tensor products fail to preserve unconditionality. Later, P\'erez-Garc\'ia and Villanueva \cite[Theorem 1.1, pp. 267-272]{Perez.Villanueva.2004} showed in 2004 that none of Grothendieck's 14 tensor norms \cite[p. 106]{Diestel.2008} satisfies the unconditionality-preserving property.
\subsection{Tests of Unconditionality}
Several tests exist to determine whether a tensor norm preserves the unconditionality of Schauder bases.
\begin{enumerate}[label=\thesubsection.\arabic*, leftmargin=4em]
\item
The First Test (Kwapie\'{n} and Pe\l czy\'{n}ski \cite[Theorem 2.2, p. 50]{Kwapien.Pelski.1970}, $1970$):\\ 
Let $\alpha$ be a tensor norm on $\ell_{2}\otimes \ell_{2}$. If $\ell_{2}\otimes_{\alpha}\ell_{2}$ admits an unconditional Schauder basis, then $\alpha$ must be equivalent to the Hilbert-Schmidt norm. For any tensor 
$u=\sum_{i,j=1}^{k,k^{'}}\lambda_{i,j} e_{i}\otimes e_{j}$ in $\ell_{2}\otimes \ell_{2}$, this norm is defined as:
$$\left\Vert u\right\Vert_{HS}
=
\left(\sum_{i,j=1}^{k,k^{'}}\vert\lambda_{i,j}\vert^{2}\right)^{1/2},
$$
where $\left(e_{m}\right)_{m\in\mathbb{N}^{*}}\subset\ell_{2}$ is a fixed orthonormal basis for $\ell_{2}$.
\item
The Second Test (P\'erez-Garc\'ia and Villanueva \cite{Perez.Villanueva.2004}, 2004):\\
If a tensor norm $\alpha$ preserves unconditionality for all pairs of Banach spaces $E$ and $F$, then $\alpha$ must coincide (up to constants) with:
\begin{enumerate}
\item
the injective norm $\varepsilon$ on $c_{0}\otimes c_{0}$ \cite[Proposition 2.3, p.270]{Perez.Villanueva.2004}, and
\item
the projective norm $\pi$ on $\ell_{1}\otimes \ell_{1}$  \cite[Proposition 2.6, p. 271]{Perez.Villanueva.2004}.
\end{enumerate}
\item
The Third Test (Diestel \cite[Proposition 4.5.8, p. 176]{Diestel.2008}, 2008):\\
Let $E$ and $F$ be real Banach spaces with unconditional bases 
$(a_{m})_{m\in\mathbb{N}^{*}}$ and 
$( x_{n})_{n\in\mathbb{N}^{*}}$
 respectively. For any tensor norm $\alpha$,  the sequence 
$\left (a_{m}\otimes x_{n}\right )_{m,n\in\mathbb{N}^{*}}$ forms an unconditional basis for $E\otimes_{\alpha}F$ if and only if $E\otimes_{\alpha}F$ is isomorphic to a Banach lattice. 
\begin{quote}
Recall on Banach Lattices: A Banach lattice is a real Banach space \(L\) equipped with a partial order \(\preceq\) satisfying the following conditions:
\begin{enumerate}
\item
Vector Space Structure Compatibility:\begin{itemize}
\item
If \(x \preceq y\), then \(x + z \preceq y + z\) for all \(z \in L\).
\item
If \(x \preceq y\), then \(\lambda x \preceq \lambda y\) for all scalar \(\lambda \geq 0\).
\end{itemize}
\item
Lattice Structure:
 For any \(x, y \in L\), there exists a least upper bound \(x \vee y\) and a greatest lower bound \(x \wedge y\).
\item 
Norm Compatibility:
If \(|x| \preceq |y|\), then \(\|x\|_{L} \leq\|y\|_{L}\), where the absolute value is defined by \(|x| = x \vee (-x)\)
    \end{enumerate}
    \end{quote}
\end{enumerate}
\subsection{Definitions}
We recall the following fundamental concepts for a Banach space $E$.
\subsubsection{Schauder Basis (SB)}
A sequence
$\left(a_{m}\right)_{,\in \mathbb{N}^{*}}\subset E$ is called a Schauder basis (SB) for $E$ if for every $x\in E$, there exists a unique scalar sequence $\left(\lambda _{m}\right)_{m\in\mathbb{N}^{*}}$
such that:
$$x=\sum_{m=1}^{\infty}\lambda_{m}a_{m}.$$
If the series $\sum_{m}\lambda_{m=1}^{\infty}a_{m}$ converges unconditionally to $x$ for every $x\in E$, then 
$\left(a_{m}\right)_{m\in \mathbb{N}^{*}}$ is called an unconditional Schauder basis (USB) for $E$.
\begin{remark}  
Every Schauder basis $\left(a_{m}\right)_{m\in \mathbb{N}^{*}}$ for $E$ naturally induces a unique biorthogonal system \(((a_m, b_{m}^{*}))_{m \in \mathbb{N}^*} \subseteq E \times E^*\), where the coefficient functionals \(b_{m}^{*}\in E^*\) satisfy \(b_{i}^{*}(a_j) = \delta_{i,j}\) for all \(i, j\in \mathbb{N}^*\). These functionals are continuous and allow the expansion of any \(x \in E\) as:
\[
x = \sum_{m=1}^{\infty} b_{m}^{*}(x)a_{m}.
\]
\end{remark}
\subsubsection{Schauder Frame (SF)}
A sequence  
$\left( (a_{m},b_{m}^{*})\right)_{m\in\mathbb{N}^{*}}\subset E\times E^{*}$
 is called a Schauder frame (SF) for $E$ (see \cite{Casazza.D.S.Z.2008}) if for every
 $x\in E$:
 \[ x=\underset{m=1}{\overset{\infty}{\sum }}
 b_{m}^{*}\left(x\right)a_{m}.\]
 If the series $\sum_{m=1}^{\infty}
 b_{m}^{*}\left(x\right)a_{m}$ converges unconditionally to $x$ for every $x\in E$, then $\left((a_{m},b_{m}^{*})\right)_{m\in\mathbb{N}^{*}}$ is called an unconditional Schauder frame (USF) for $E$.
\subsubsection{Besselian Schauder Frame (BSF)}
 A sequence 
$\left( (a_{m},b_{m}^{*})\right)_{n\in\mathbb{N}^{*}}\subset E\times E^{*}$ is called a Besselian Schauder frame (BSF) for $E$ (see \cite{karkri.zoubeir.2021}) if:
\begin{enumerate}
\item
It is a Schauder frame (SF) for $E$, and
\item
There exists a constant $B>0$ such that for every $x\in E$ and $f^{*}\in E^{*}$:
\begin{equation*}
\underset{m=1}{\overset{\infty}{\sum }}\left\vert b_{m}^{*}\left(
x\right)\right\vert\left\vert f^{*}\left(a_{m}\right)\right\vert
\leq B\left\Vert x\right\Vert_{E}\left\Vert f^{*}\right\Vert
_{E^{*}}.
\end{equation*}
\end{enumerate}
\begin{remark} Let $\mathcal{F}:=\left( \left(
a_{m},b_{m}^{*}\right)\right)_{m\in\mathbb{N}^{*}}$ be a BSF for $E$, the quantity 
\begin{equation*}
\mathcal{L}_{\mathcal{F}}
:=\sup_{\substack{\Vert x\Vert_{E} \leq 1 \\ \Vert f^{*}\Vert_{E^{*}}\leq 1}}
\sum_{m=1}^{\infty}
\left\vert b_{m}^{*}\left(x\right)\right\vert
 \left\vert f^{*}\left(a_{m}\right)\right\vert
\end{equation*}
is finite. Moreover, for every $\left(x,f^{*}\right)\in E\times E^{*}$, the
following inequality holds:
\begin{equation*}
\sum_{m=1}^{\infty}
\left\vert b_{m}^{*}\left(x\right)
\right\vert\left\vert f^{*}\left(a_{m}\right)\right\vert\leq 
\mathcal{L}_{\mathcal{F}}\left\Vert x\right\Vert_{E}
\left\Vert f^{*}\right\Vert_{E^{*}}.
\end{equation*}
The constant $\mathcal{L}_{\mathcal{F}}$ is called the Besselian constant of the
frame $\mathcal{F}$.
\end{remark}

In this paper, for a given pair of Banach spaces $\left(E,\Vert\cdot\Vert_{E}\right )$
and $\left(F,\Vert\cdot\Vert_{F}\right)$ with Besselian Schauder frames $\mathcal{F}=\left((a_{m},b_{m}^{*})\right)_{m\in \mathbb{N}^{*}}$ and $\mathcal{G}=\left((x_{n},y_{n}^{*})\right)_{n\in \mathbb{N}^{*}}$, we proceed in two steps: 
\begin{enumerate}
\item
First, we construct equivalent norms
$\Vert\cdot\Vert^{Bess}_{E,\mathcal{F}}$ on $E$ and $\Vert\cdot\Vert^{Bess}_{F,\mathcal{G}}$ on $F$.
\item
Using these new norms, we define a reasonable crossnorm $\alpha_{\mathcal{F},\mathcal{G}}^{Bess}$ on $E\otimes F$ such that:
\begin{enumerate}
\item
If $\mathcal{F}$ and $\mathcal{G}$ are unconditional Schauder bases, then:
\begin{itemize}
\item
The tensor product sequence 
$\left((a_{m}\otimes x_{n},b_{m}^{*}\otimes y_{n}^{*})\right)_{m,n\in \mathbb{N}^{*}},$
ordered via the square ordering, forms a USB for $E\otimes_{\alpha_{\mathcal{F},\mathcal{G}}^{Bess}}F$.
\item
If $E$ and $F$ are real spaces, then $\left(E,\Vert\cdot\Vert^{Bess}_{E,\mathcal{F}}\right)$, $\left(F,\Vert\cdot\Vert^{Bess}_{F,\mathcal{G}}\right)$, and $\left(E\otimes_{\alpha_{\mathcal{F},\mathcal{G}}^{Bess}}F,\alpha_{\mathcal{F},\mathcal{G}}^{Bess}\right)$ are Banach lattices.
\end{itemize}
\item
If  $\mathcal{F}$ and $\mathcal{G}$ are the canonical USB of $c_{0}$ or $\ell_{1}$. Then
\begin{itemize}
\item
$\alpha_{\mathcal{F},\mathcal{G}}^{Bess}$ conicides with the projective norm $(\pi)$ on $\ell_{1}\otimes \ell_{1}$. 
\item
$\alpha_{\mathcal{F},\mathcal{G}}^{Bess}$ conicides with the injective norm $(\varepsilon)$ on $c_{0}\otimes c_{0}$. 
\end{itemize}
\end{enumerate}
\end{enumerate}
This new norm $\alpha_{\mathcal{F},\mathcal{G}}^{Bess}$ differs fundamentally from the classical Chevet-Saphar norms $d_{p},g_{p}$ (see \cite[p.133]{Raymond.A.Ryan}) and  the norms defined by S. Kwapie\'{n} and A. Pe\l czy\'{n}ski in \cite[Definition 1.2, p. 47]{Kwapien.Pelski.1970} in several key aspects:
\begin{itemize}
\item 
The Chevet-Saphar and Kwapie\'{n}-Pe\l czy\'{n}ski norms are defined using the structure of $\ell_p$-spaces and are applicable to general tensor products of Banach spaces.
\item
In contrast, \(\alpha_{\mathcal{F},\mathcal{G}}^{\text{Bess}}\) is explicitly constructed for pairs of Banach spaces \(E\) and \(F\) equipped with specific BSFs (or USBs) \(\mathcal{F}\) and \(\mathcal{G}\). It leverages the  structure of these frames (via the Besselian norms $\Vert\cdot\Vert^{Bess}_{E,\mathcal{F}}$ and $\Vert\cdot\Vert^{Bess}_{F,\mathcal{G}}$) to define a norm intrinsically tailored to the geometry of \(E\) and \(F\).
\item
This approach allows \(\alpha_{\mathcal{F},\mathcal{G}}^{\text{Bess}}\)---although this norm is generally not uniform---to capture finer properties of the spaces, such as unconditionality and lattice structure, which are not generally accessible via  the classical norms.
\end{itemize} 
\section{Notations and Preliminaries}
\label{Notations.preliminaries}
Throughout this paper, $\left(E,\Vert\cdot\Vert_{E}\right)$ (resp. $\left(F,\Vert\cdot\Vert_{F}\right)$) denotes a
 Banach space (real or complex), and $E^{*}$ (resp. $F^{*}$) represents its topological dual.
For $p \in (1, +\infty)$, the conjugate exponent is defined as:
$$p^{*}=\dfrac{p}{p-1}.$$
\subsection{Sequences and Dual Sequences:}
\begin{enumerate}[label=\thesubsection.\arabic*, leftmargin=4em]
\item
Let  $\mathcal{F}:=\left( \left(a_{m},b_{m}^{*}\right)\right)_{m\in \mathbb{N}^{*}}$ be a fixed sequence in $E\times E^{*}$.
\item
Similarly, let $\mathcal{G}:=\left( \left( x_{n},y_{n}^{*}\right) \right)_{n\in \mathbb{N}^{*}}$ be a fixed sequence in $F\times F^{*}$.
 \item
Let $G$ be a Banach space and  $\left((z_{k},h_{k}^{*})\right)_{k\in\mathbb{N}^{*}}$ be a sequence in $G\times G^{*}$. The sequence $\left((z_{k},h_{k}^{*})\right)_{k\in\mathbb{N}^{*}}$ is called a biorthogonal system if 
 $$h_{i}^{*}(z_{j})=\delta_{i,j},\;\;(i,j\in\mathbb{N}^{*}),$$
where $\delta_{i,j}$ denotes the Kronecker delta.
 \item
Denote by $\mathcal{S}$ the set of all sign sequences:  
\[
\mathcal{S} = \left\{s=\left (s_{m}\right )_{m \in\mathbb{N}^{*}}\subset\mathbb{K}: \vert s_{m}\vert=1\right\}.
\]  
We equip $\mathcal{S}$ with the metric:  
\[
d(s,t)=\sum_{m=1}^\infty \frac{|s_m - t_m|}{2^m}.
\]  
For each component space $\{z \in \mathbb{K}: \vert z\vert= 1\}$ is compact (closed and bounded in $\mathbb{K}$). By Tychonoff's theorem, the product of countably many compact sets is compact in the product topology. The metric $d$ induces this topology since $s^{(k)} \to s$ in $(\mathcal{S},d)$ if and only if $s_{m}^{(k)} \to s_{m}$ in $\mathbb{K}$ for all $m$. Thus, $\mathcal{S}$ is compact.  
\end{enumerate}
\subsection{Classical Spaces}
\begin{enumerate}[label=\thesubsection.\arabic*, leftmargin=4em]
\item
Let $\mathbb{K}=\mathbb{R}$ or $\mathbb{K}=\mathbb{C}$. We denote by $c_{0}$ the Banach space of all sequences in $\mathbb{K}$ converging to zero:
$$c_{0}:=\left\lbrace \left(\lambda_{i}\right)_{i\in\mathbb{N}^{*}}\subset \mathbb{K}:
\lim_{i\to \infty}\lambda_{i}=0
\right\rbrace ,$$
equipped with the supremum norm:
$$
\left\Vert\left(\lambda_{i}\right)_{i\in\mathbb{N}^{*}}\right\Vert_{c_{0}}
=
\underset{i\in\mathbb{N}^{*}}{\sup} \vert\lambda_{i}\vert .
$$
Its canonical Schauder basis is denoted by $\left(e_{n}\right)_{n\in \mathbb{N}^{*}}$, where $e_{n}=\left(\delta_{n,k}\right)_{k\in \mathbb{N}^{*}}$,
 with biorthogonal functionals $\left(e_{n}^{*}\right)_{n\in \mathbb{N}^{*}}$, \item
We denote by $c_{0}(E)$ the Banach space of all $E$-valued sequences converging to zero:
$$c_{0}(E):=\left\lbrace \left(x^{i}\right)_{i\in\mathbb{N}^{*}}\subset E:
\lim_{i\to \infty}\Vert x^{i}\Vert_{E}=0
\right\rbrace ,$$
equipped with the supremum norm:
$$
\left\Vert\left(x^{i}\right)_{i\in\mathbb{N}^{*}}\right\Vert_{c_{0}(E)}
=
\underset{i\in\mathbb{N}^{*}}{\sup} \Vert x^{i}\Vert_{E}.
$$
\item
We denote by $\ell_{\infty}$ the Banach space of all bounded sequences in $\mathbb{K}$:
$$\ell_{\infty}:=\left\lbrace\left(\lambda_{i}\right)_{i\in\mathbb{N}^{*}}\subset \mathbb{K}:
\sup_{i\in\mathbb{N}^{*}}\vert\lambda_{i}\vert<\infty
\right\rbrace ,$$
equipped with the supremum norm:
$$
\left\Vert\left(\lambda_{i}\right)_{i\in\mathbb{N}^{*}}\right\Vert_{\ell_{\infty}}
=
\underset{i\in\mathbb{N}^{*}}{\sup} \vert\lambda_{i}\vert .
$$
\item
We denote by $\ell_{p}$ (for $1\leq p<\infty$) the Banach space of all absolutely p-summable sequences in $\mathbb{K}$:
$$\ell_{p}:=\left \lbrace \left(\lambda_{i}\right) _{i\in\mathbb{N}^{*}}\subset\mathbb{K}:
\sum_{i=1}^{\infty }\left \vert\lambda_{i}\right\vert^{p}
<\infty
\right  \rbrace ,$$
equipped with the norm:
$$
\left\Vert(\lambda_{i})_{i\in\mathbb{N}^{*}}\right\Vert_{\ell_{p}}
=
\left(\sum_{i=1}^{\infty }\left \vert\lambda_{i}\right \vert^{p}
\right )^{1/p}.$$
Its canonical Schauder basis is denoted by $\left(e_{n}\right)_{n\in \mathbb{N}^{*}}$, with biorthogonal functionals $\left(e_{n}^{*}\right)_{n\in \mathbb{N}^{*}}$.
\item We denote by $\ell_{p}(E)$ the Banach space of strongly p-summable sequences in $E$:
$$\ell_{p}(E):=\left \lbrace \left( x^{i}\right) _{i\in \mathbb{N}^{* }}\subset E:
\sum_{i=1}^{\infty }\left \Vert x^{i}\right\Vert_{E}^{p}
<\infty
\right  \rbrace ,$$
equipped with the norm:
$$
\left\Vert(x^{i})_{i\in\mathbb{N}^{*}}\right\Vert_{\ell_{p}(E)}
=
\left(\sum_{i=1}^{\infty }\left \Vert x^{i}\right\Vert_{E}^{p}
\right )^{1/p}.$$
\item
We denote by $L\left( E\right) $ the Banach space of all bounded linear
operators $T:E\rightarrow E$, equipped with the operator norm:
\[
\|T\|_{L(E)}:= \sup_{\substack{\|x\|_{E}\leq 1}} \|T(x)\|_{E}.
\]
\end{enumerate}
\subsection{Tensor Product of Banach Spaces (see, \cite{Diestel.2008}, \cite{Raymond.A.Ryan} and \cite{R.Schatten.1950}).}
\begin{enumerate}[label=\thesubsection.\arabic*, leftmargin=4em]
\item
The algebraic tensor product of the spaces $E$ and $F$, denoted $E\otimes F$, is the linear span of all elementary tensors:
$$E\otimes F:=span\left\{ x\otimes y: x\in E, y\in F\right\}.$$
Each elementary tensor $x\otimes y$ 
defines a bounded linear operator:
\begin{align*}
x\otimes y:  E^{*} \longrightarrow  F,\;\;
f^{*} \longmapsto f^{*}(x)y.
\end{align*} 
\item
The algebraic tensor product of dual spaces $E^{*}\otimes F^{*}$ is the linear span:
$$E^{*}\otimes F^{*}:=span\left \{f^{*}\otimes g^{*}: f^{*}\in E^{*}, g^{*}\in F^{*}\right\}.$$
Each elementary tensor $f^{*}\otimes g^{*}$   induces a linear functional on $E\otimes F$:
  \begin{align*}
f^{*}\otimes g^{*}:\;  E\otimes F \longrightarrow   \mathbb{K},\;\;
\sum_{r=1}^{R} x^{r}\otimes y^{r}\longmapsto\sum_{r=1}^{R}f^{*}(x^{r})g^{*}(y^{r}).
\end{align*}
\item
 A norm $\alpha$ on $E\otimes F$  is called a reasonable crossnorm   if it satisfies the following two conditions:
 \begin{enumerate}
 \item
 For every $x\in E$ and $y\in F$,
 $$\alpha (x\otimes y)\leq\Vert x\Vert_{E}\Vert y\Vert_{F}.$$
 \item
For every $f^{*}\in E^{*}$, $g^{*}\in F^{*}$, and every $u\in E\otimes F$,
\begin{align*}
\left\vert f^{*}\otimes g^{*}\left(u\right)\right\vert
     &\leq 
     \Vert f^{*}\Vert_{E^{*}}\Vert g^{*}\Vert _{F^{*}}\alpha\left(u\right).
\end{align*}
\end{enumerate}
(see \cite[Definition 1.1.1, p. 5]{Diestel.2008} or \cite[p. 127]{Raymond.A.Ryan}).
\item
A reasonable crossnorm $\alpha$ on $E\otimes F$ is called uniform if for all bounded linear operators $S\in L(E)$ and $T\in L(F)$, the linear operator: 
\begin{align*}
S\otimes T:\;  E\otimes F \longrightarrow E\otimes F ,\;\;
\sum_{r=1}^{R} x^{r}\otimes y^{r}\longmapsto\sum_{r=1}^{R}S(x^{r})\otimes T(y^{r})
\end{align*}
is bounded with respect to $\alpha$ and its norm satisfies:
\begin{align*}
\left\Vert S\otimes T\right\Vert_{L(E\otimes F)}
     &\leq 
     \Vert S\Vert _{L(F)}\Vert T\Vert_{L(E)}.
\end{align*}
Equivalently, for every tensor $u\in E\otimes F$,
\begin{align*}
\alpha\left(S\otimes T(u)\right)
     &\leq 
     \Vert S\Vert _{L(F)}\Vert T\Vert_{L(E)}
 \alpha\left(u\right).
\end{align*}
 \item 
The completion of $E\otimes F$ with 
 respect to $\alpha$ is denoted by $E\otimes_{\alpha}F$. 
 \item 
The injective tensor norm $\varepsilon$ is defined for a tensor $u=\sum_{r=1}^{R}x^{r}\otimes y^{r}$ in $E\otimes F$ by 
\[\varepsilon\left(u\right)
=\sup_{\substack{\|f^*\|_{E^*} \leq 1 \\ \|g^*\|_{F^*} \leq 1}}
\left \vert \underset{r=1}{\overset{R}{\sum}}f^{*}(x^{r})g^{*}(y^{r})\right\vert .
\]
The completion of $E\otimes F$ with 
 respect to $\varepsilon$ is denoted by $E\otimes_{\varepsilon}F$. 
 \item
The projective tensor norm $\pi$ is defined for a tensor $u\in E\otimes F$ by 
\[ \pi(u)=\inf\left \{\underset{r=1}{\overset{R}{\sum}}\Vert x^{r}\Vert_{E}\Vert y^{r}\Vert_{F} :
u=\underset{r=1}{\overset{R}{\sum }}x^{r}\otimes y^{r}\right\},\]
where the infimum is taken over all representations of $u$ as a finite sum of elementary tensors. The completion of $E\otimes F$ with 
 respect to $\pi$ is denoted by $E\otimes_{\pi}F$. 
\end{enumerate}
\section{Main Results}
\begin{lemma}
 Let \(\sum_{m=1}^\infty z_m\) be a series in a Banach space \(E\). The following statements are equivalent:
 \begin{enumerate}
 \item
The series converges unconditionally.
\item
For every \(\varepsilon > 0\), there exists \(M \in \mathbb{N}^{*}\) such that for all finite sets \(A \subseteq \{M+1, M+2, \ldots\}\) and all scalars \(\lambda_m \in \mathbb{K}\) with \(|\lambda_m| \leq 1\):
   \[
   \left\| \sum_{m \in A} \lambda_m z_m \right\|_E < \varepsilon.
   \]
 \end{enumerate}
\end{lemma}
\textbf{Proof.} 
The implication $(2)\Rightarrow (1)$ is classical. Indeed, if condition (2) holds, then in particular (by taking $\lambda_{m}=1$ for all $m$) the series satisfies the Cauchy criterion and is therefore convergent. The unconditional convergence then follows from a standard characterization \cite[Theorem 16.1, p. 461]{I.Singer.I}.\\
 Now, assume that the series converges unconditionally.  Suppose, for contradiction, that condition (2) fails. Then there exists \(\varepsilon_0 > 0\) such that for every \(M\in \mathbb{N}^{*}\), there exist a finite set \(A_M \subseteq \{M+1, M+2, \ldots\}\) and scalars \((\lambda_m^M)_{m\in A_{M}}\) with \(|\lambda_m^M| \leq 1\) satisfying:
\[
\left\| \sum_{m \in A_M} \lambda_m^M z_m \right\|_E \geq \varepsilon_0.
\]

We  now construct a bounded sequence $(\lambda_m)_{m\in\mathbb{N}^{*}}$ such that the series \(\sum_{m=1}^{\infty}\lambda_m z_m\) diverges, contradicting the classical characterization of unconditional convergence.
\begin{itemize}
\item
Let \(M_1 = 1\). Choose a finite set \(A_{1} \subseteq \{M_1+1, M_1+2, \ldots\}\) and scalars \((\lambda_m^{1})_{m\in A_{1}}\) with \(|\lambda_m^{1}| \leq 1\) such that
$$\left\| \sum_{m \in A_{1}} \lambda_m^{1} z_m \right\|_E \geq \varepsilon_0.$$
\item
Set \(M_2 = \max A_{1} + 1\). Choose a finite set \(A_{2} \subseteq \{M_2+1, M_2+2, \ldots\}\) and scalars \((\lambda_m^{2})_{m\in A_{2}}\) with \(|\lambda_m^{2}| \leq 1\) such that
$$\left\| \sum_{m \in A_{2}} \lambda_m^{2} z_m \right\|_E \geq \varepsilon_0.$$
\item
Proceed inductively: having defined \(M_k\) and \(A_k\), set \(M_{k+1} = \max A_{k} + 1\), and choose a finite set \(A_{k+1} \subseteq \{M_{k+1}+1, M_{k+1}+2,\ldots\}\) and scalars \((\lambda_m^{k+1})_{m\in A_{k+1}}\) with \(|\lambda_m^{k+1}| \leq 1\) such that
$$\left\| \sum_{m \in A_{k+1}} \lambda_m^{k+1} z_m \right\|_E \geq \varepsilon_0.$$
 \end{itemize}
Define the sequence \((\lambda_m)_{m\in\mathbb{N}^{*}}\) by:
\[
\lambda_m = 
\begin{cases}
\lambda_m^{k} & \text{if } m \in A_{k} \text{ for some } k, \\
0 & \text{otherwise}.
\end{cases}
\]
By construction, \(|\lambda_m| \leq 1\) for all $m$, so \((\lambda_m)_{m\in\mathbb{N}^{*}}\) is bounded. Now consider the partial sums over the finite sets \(S_k = \bigcup_{j=1}^k A_{j}\). The sets $(A_k)_{k\in\mathbb{N}^{*}}$ are finite, disjoint, and their union is increasing. For each $k\geq 2$, we have:
\[
\left\| \sum_{m \in S_k} \lambda_m z_m - \sum_{m \in S_{k-1}} \lambda_m z_m \right\|_E = \left\| \sum_{m \in A_{k}} \lambda_m^{k} z_m \right\|_E \geq \varepsilon_0.
\]
This shows that the sequence of partial sums $\left(\sum_{m\in S_{k}} \lambda_m z_m\right)_{k\in\mathbb{N}^{*}}$ is not Cauchy. Therefore, the series $\sum \lambda_m z_m$ diverges, contradicting the assumption of unconditional convergence.\qed\\

The following theorem is inspired by  \cite[p. 32-33]{A.Abramovich.2002}.
\begin{theorem}
\label{unc.impl.bess}
Suppose that $\mathcal{F}$ is an USF for $E$. Then $\mathcal{F}$ is a BSF for $E$.
\end{theorem}
\textbf{Proof}.
For each  $M\in \mathbb{N}^{*}$  and 
 $(s_{m})_{m\in\mathbb{N}^{*}}\in\mathcal{S}$, the linear operator
  $$T_{M,s}:E\rightarrow E,\; x\mapsto\sum_{m=1}^{M}s_{m}b_{m}^{*}(x)a_{m}$$ 
is  bounded. By the Banach-Steinhaus theorem, the linear operator
  $$T_{s}:E\rightarrow E,\; x\mapsto\sum_{m=1}^{\infty} 
  s_{m}b_{m}^{*}(x)a_{m}$$ is bounded.
We now show that for each $x\in E$, the mapping 
 $s\mapsto T_{s}(x) $ is continuous on the compact metric space
 $\mathcal{S}$.\\
  Fix $x\in E$ and $\varepsilon > 0$. By unconditionality of $\mathcal{F}$, there exists $M\in\mathbb{N}^{*}$ such that for all finite sets 
 $A \subseteq\{M+1, M+2, \dots\}$ and all scalars $\lambda_{m}\in\mathbb{K}$ with 
 $|\lambda_m| \leq 1$:  
  \[
  \left\|\sum_{m\in A}\lambda_m b_m^*(x) a_m \right\|_E < \frac{\varepsilon}{4}.
  \]  
Now, for \(s, t \in \mathcal{S}\), we estimate the difference:
  \[
  \|T_s(x) - T_t(x)\|_E \leq \underbrace{\left\| \sum_{m=1}^M(s_m - t_m) b_m^*(x) a_m \right\|_E}_{\text{(I)}} + \underbrace{\left\| \sum_{m=M+1}^\infty (s_m - t_m) b_m^*(x) a_m \right\|_E}_{\text{(II)}}.
  \]  
Estimate (II): Since \(|s_m - t_m| \leq 2\), define \(\lambda_m = \frac{s_m - t_m}{2}\) (\(|\lambda_m| \leq 1\)). Then,
    \[
    \text{(II)} \leq 2  \sup_{A \subseteq \{M+1,\dots\}} \left\| \sum_{m \in A} \lambda_m b_m^*(x) a_m \right\|_E \leq 2  \frac{\varepsilon}{4} = \frac{\varepsilon}{2}.
    \]  
Estimate (I): Let \(C = \max_{1 \leq m \leq M} \|b_m^*(x) a_m\|_E\). Choose 
   $$\delta = \frac{\varepsilon}{2^{M+1}MC}.$$
If \(d(s,t) < \delta\), then for each \(m \leq M\),
$$|s_m - t_m|<2^m d(s,t)< 2^m \delta \leq 2^M\delta .$$
Hence,
  \[
    \text{(I)} \leq \sum_{m=1}^M |s_m - t_m|  C <M(2^M\delta) C = \frac{\varepsilon}{2}.
    \]  
  Combining (I) and (II), we conclude that:
$$\|T_s(x) - T_t(x)\|_E < \varepsilon$$
 whenever \(d(s,t) < \delta\). Thus, $s\mapsto T_{s}(x) $ is continuous on 
 $\mathcal{S}$.\\
  Since $\mathcal{S}$ is compact, the set
 $\left\lbrace \Vert T_{s}(x)\Vert_{E}: s\in \mathcal{S}\right\rbrace$ is bounded. 
By the uniform boundedness principle,  there is a constant $\mathcal{D}>0$ such that
$$\left \Vert \underset{m=1}{\overset{\infty}{\sum}}
s_{m}b_{m}^{*}(x)a_{m}\right\Vert_{E}
\leq
\mathcal{D}\Vert x\Vert_{E},
\;\;\;\left (x\in E,\;(s_{m})_{m\in\mathbb{N}^{*}}\in \mathcal{S}\right ).$$
Now fix $x\in E$ and $f^{*}\in E^{*}$ with $\Vert f^{*}\Vert_{E^{*}}\leq 1$. Choose a sequence of signs
$s=\left( s_{m}\right)_{m\in \mathbb{N}^{*}}$ such that 
$$\vert b_{m}^{*}(x)f^{*}(a_{m})\vert=s_{m}b_{m}^{*}(x)f^{*}(a_{m}).$$
Then:
\begin{align*}
\underset{m=1}{\overset{\infty}{\sum}}
\left\vert b_{m}^{*}(x)f^{*}(a_{m})\right\vert 
& =\overset{\infty}{\underset{m=1}{\sum }}
s_{m}b_{m}^{*}(x)f^{*}(a_{m})
=\left\vert f^{*}\left(\overset{\infty }{\underset{m=1}{\sum}}
s_{m}b_{m}^{*}(x)a_{m}\right)\right\vert\\
&=\left\vert f^{*}\left(T_{s}(x)\right)\right\vert 
\leq\left\Vert T_{s}(x)\right\Vert_{E}
\leq\left\Vert T_{s}\right\Vert_{L(E)}\cdot\left\Vert x\right\Vert_{E}\\
&\leq \mathcal{D}\left\Vert x\right\Vert_{E}.
\end{align*}
For general $f^{*}\in E^{*}$, scale by $\Vert f^{*}\Vert_{E^{*}}$ to obtain:
$$
\underset{m=1}{\overset{\infty}{\sum}}
\left\vert b_{m}^{*}(x)f^{*}(a_{m})\right\vert 
\leq \mathcal{D}\left \Vert  x\right\Vert_{E}\left\Vert f^{*}\right \Vert_{E^{*}}, 
\;\;\; (x\in E, f^{*}\in E^{*}).
$$
Thus $\mathcal{F}$ is a Besselian Schauder frame for $E$.\qed\\
 \begin{lemma} 
 \label{norm.dual}
 Suppose that $\mathcal{F}$ is a BSF for $E$. Then
 \begin{enumerate}
 \item
The mapping:
$$\Vert\cdot\Vert^{Bess}_{E,\mathcal{F}} :E\rightarrow \mathbb{R}^{+},\;\;x\mapsto \underset{\Vert f^{*}\Vert_{E^{*}}\leq 1}{\sup}\underset{m=1}{\overset{\infty}{\sum}}
 \left\vert b_{m}^{*}(x)f^{*}(a_{m})\right\vert,
   $$
defines a norm (Besselian norm) on $E$ equivalent to the original norm $\left\Vert x\right\Vert_{E}$, satisfying
   $$\left\Vert x\right\Vert_{E}
   \leq \Vert x\Vert^{Bess}_{E,\mathcal{F}}
   \leq \mathcal{L}_{\mathcal{F}}\left \Vert x\right\Vert_{E} ,\;\;(x\in E),$$
where $\mathcal{L}_{\mathcal{F}}$ is the Besselian constant of the frame $\mathcal{F}$.
\item
The dual norm on $E^{*}$ induced by $\Vert\cdot\Vert_{E,\mathcal{F}}^{Bess}$ is given by
$$\Vert f^{*}\Vert_{E^{*},\mathcal{F}}^{Bess}:=\underset{\Vert x\Vert^{Bess}_{E,\mathcal{F}}\leq 1}{\sup}
 \left \vert f^{*}(x)\right\vert,\;\;(f^{*}\in E^{*}),
   $$ 
and satisfies
\begin{align}
\label{norms.dual.equi}
\mathcal{L}_{\mathcal{F}}^{-1}\left \Vert f^{*}\right\Vert_{E^{*}}
   \leq \Vert f^{*}\Vert_{E^{*},\mathcal{F}}^{Bess}
   \leq \left \Vert f^{*}\right\Vert_{E^{*}} ,\;\;(f^{*}\in E^{*}).
\end{align}
\end{enumerate}  
\end{lemma}
\textbf{Proof.}
\begin{enumerate}
\item
For each $x\in E$, the Schauder basis decomposition $x=\sum_{m=1}^{\infty}b_{m}^{*}(x)a_{m}$, yields:
\begin{align*}
\Vert x\Vert_{E}
&=
\underset{\Vert f^{*}\Vert_{E^{*}}\leq 1}{\sup}
 \left\vert f^{*}(x)\right\vert
 =
\underset{\Vert f^{*}\Vert_{E^{*}}\leq 1}{\sup}
 \left \vert
 \sum_{m=1}^{\infty} b_{m}^{*}(x)f^{*}(a_{m})
 \right\vert\\ 
 &\leq
\underset{\Vert f^{*}\Vert_{E^{*}}\leq 1}{\sup}
 \sum_{m=1}^{\infty} \left \vert
 b_{m}^{*}(x)f^{*}(a_{m})
 \right\vert 
\stackrel{\scriptstyle\text{def}}{=}
  \Vert x\Vert^{Bess}_{E,\mathcal{F}}
 \leq
 \mathcal{L}_{\mathcal{F}}
\Vert x\Vert_{E}.
\end{align*}
This proves the first item.
\item
For each $f^{*}\in E^{*}$, we have:
\begin{align*}
\mathcal{L}_{\mathcal{F}}^{-1}\left\Vert f^{*}\right\Vert_{E^{*}}
&=
\mathcal{L}_{\mathcal{F}}^{-1}
\underset{\Vert x\Vert_{E}\leq 1}{\sup}
 \left \vert f^{*}(x)\right\vert
 =
\underset{\mathcal{L}_{\mathcal{F}}\Vert x\Vert_{E}\leq 1}{\sup}
 \left \vert f^{*}(x)\right\vert\\ 
 &\leq
\underset{\Vert x\Vert^{Bess}_{E,\mathcal{F}}\leq 1}{\sup}
 \left \vert f^{*}(x)\right\vert
\stackrel{\scriptstyle\text{def}}{=}
 \Vert f^{*}\Vert_{E^{*},\mathcal{F}}^{Bess}
 \leq
 \underset{\Vert x\Vert_{E}\leq 1}{\sup}
 \left\vert f^{*}(x)\right\vert
 =\Vert f^{*}\Vert_{E^{*}}.
\end{align*}
This proves the second item.
\end{enumerate} \qed
\begin{example}
\label{example.norm.bess=norm}
Suppose that $E$ is $\ell_{p}$ (for $1\leq p <\infty$) or $E=c_{0}$, and let $\mathcal{F}$ be the canonical USB of $E$. Then,
\begin{align}
\label{norm.bess=norm}
\Vert r\Vert_{E,\mathcal{F}}^{Bess}
=
\left\Vert r\right\Vert_{E} , \;\;\;\;(r\in E).
\end{align}
Consequently, the Besselian constant is 
$\mathcal{L}_{\mathcal{F}}=1$.
\end{example}
\textbf{Proof.}
 Let $r=(r_{m})_{m\in\mathbb{N}^{*}}\in E$. We have
\begin{align*}
\left\Vert r\right\Vert_{E}
&=
\underset{\Vert f^{*}\Vert_{E^{*}}\leq 1}{\sup}
\left\vert f^{*}(r)\right\vert
=
\underset{\Vert f^{*}\Vert_{E^{*}}\leq 1}{\sup}
\left\vert
\underset{m=1}{\overset{\infty}{\sum}}
  b_{m}^{*}(r)f^{*}(a_{m})
  \right\vert\\
&\leq
  \underset{\Vert f^{*}\Vert_{E^{*}}\leq 1}{\sup}
\underset{m=1}{\overset{\infty}{\sum}}
\left \vert
b_{m}^{*}(r)f^{*}(a_{m})
\right\vert
\stackrel{\scriptstyle\text{def}}{=}
  \Vert r\Vert_{E,\mathcal{F}}^{Bess}.
\end{align*}
On the other hand, it follows from:
\begin{itemize}
\item
H\"{o}lder's inequality if $E=\ell_{p}$, and 
\item
 the fact that 
$\vert r_{n}\vert\leq\Vert r\Vert_{E} $ for each $n\in \mathbb{N}^{*}$ if $E=c_{0}$,
\end{itemize} 
that:
\begin{align*}
\Vert r\Vert_{E,\mathcal{F}}^{Bess}
&=
\underset{\Vert f^{*}\Vert_{E^{*}}\leq 1}{\sup}
\underset{m=1}{\overset{\infty}{\sum}}
\left \vert
  b_{m}^{*}(r)f^{*}(a_{m})
  \right\vert
 =
 \underset{\Vert f^{*}\Vert_{E^{*}}\leq 1}{\sup}
\underset{m=1}{\overset{\infty }{\sum}}
   \left\vert
   r_{m}y_{m}^{*}
  \right\vert
\leq
\Vert r\Vert_{E}.
\end{align*}
Consequently,
$$\Vert r\Vert_{E,\mathcal{F}}^{Bess}
=
\left\Vert r\right\Vert_{E} , \;\;\;\;(r\in E).$$
It follows directly from the definition of 
$\mathcal{L}_{\mathcal{F}}$ and inequality \eqref{norm.bess=norm} that $\mathcal{L}_{\mathcal{F}}=1$,
which finishes the proof.\qed\\
\begin{remark}
Suppose that $\mathcal{F}$ is a USB for $E$. Then, 
$$
\left\Vert x\right\Vert_{E,\mathcal{F}}^{Bess}=\left\Vert\sum_{m=1}^{\infty}\varepsilon_{m}b_{m}^{*}(x)a_{m}\right\Vert_{E,\mathcal{F}}^{Bess}=\left\Vert\sum_{m=1}^{\infty}\vert b_{m}^{*}(x)\vert a_{m}\right\Vert_{E,\mathcal{F}}^{Bess},
$$
for every $x=\sum_{m=1}^{\infty} b_{m}^{*}(x)a_{m}$ in $E$ and every sequence of signs $(\varepsilon_{m})_{m\in\mathbb{N}^{*}}\in \mathcal{S}$.
\end{remark}
\begin{definition}
\label{defi-new-norm}
Suppose that $\mathcal{F}$ and $\mathcal{G}$ are  BSF for 
$\left(E,\Vert\cdot\Vert_{E}\right)$ and $\left (F,\Vert\cdot\Vert_{F}\right)$, respectively.
Define the mapping $\alpha_{\mathcal{F},\mathcal{G}}^{Bess}: E\otimes F\rightarrow\mathbb{R}^{+}$ by:
$$
\alpha_{\mathcal{F},\mathcal{G}}^{Bess}\left(u\right)
:=
\inf
\left \lbrace
Bess\left(\sum_{r=1}^{R} x^{r}\otimes y^{r}\right):
u=\sum_{r=1}^{R}  x^{r}\otimes y^{r}
\right \rbrace , (u\in E\otimes F),
$$
where the infimum is taken over all finite representations of $u$, and
$$
Bess\left (\sum_{r=1}^{R}x^{r}\otimes y^{r}\right)
:=
\sup_{\substack{\|f^*\|_{E^*} \leq 1 \\ \|g^*\|_{F^*} \leq 1}}
\underset{r=1}{\overset{R}{\sum}} 
\underset{m,n=1}{\overset{\infty}{\sum}} 
 \left\vert  b_{m}^{*}(x^{r})y_{n}^{*}(y^{r})f^{*}(a_{m})g^{*}\left(x_{n}\right)    \right\vert . 
$$
We refer to $\alpha_{\mathcal{F},\mathcal{G}}^{Bess}$ as the Besselian crossnorm associated with $\mathcal{F}$ and $\mathcal{G}$.
\end{definition}
\begin{proposition}
\label{New-crossnorm}
Under the hypotheses of definition \ref{defi-new-norm}, the following holds:
\begin{enumerate}
\item
For each $u\in E\otimes F$ we have:
\begin{equation}
\label{New-Norm-Between-inj-proj}
\mathcal{L}_{\mathcal{F}}^{-1}
\mathcal{L}_{\mathcal{G}}^{-1}\varepsilon_{\mathcal{F},\mathcal{G}}\left (u\right)
\leq
\alpha_{\mathcal{F},\mathcal{G}}^{Bess}\left (u\right)
\leq
\pi_{\mathcal{F},\mathcal{G}}\left(u\right),
\end{equation}
where $\varepsilon_{\mathcal{F},\mathcal{G}}$ (resp. $\pi_{\mathcal{F},\mathcal{G}}$) denotes  the injective (resp. projective) tensor norm on the algebraic tensor product of
$\left (E,\Vert\cdot\Vert^{Bess}_{E,\mathcal{F}}\right)$ 
and 
 $\left(F,\Vert\cdot\Vert^{Bess}_{F,\mathcal{G}}\right)$.
\item
The mapping $\alpha_{\mathcal{F},\mathcal{G}}^{Bess}$ defines a reasonable crossnorm on the tensor product of 
$\left(E,\Vert\cdot\Vert^{Bess}_{E,\mathcal{F}}\right)$ 
and 
 $\left(F,\Vert\cdot\Vert^{Bess}_{F,\mathcal{G}}\right)$. Consequently,
 \begin{equation}
\label{2-New-Norm-Between-inj-proj}
\varepsilon_{\mathcal{F},\mathcal{G}}\left (u\right)
\leq
\alpha_{\mathcal{F},\mathcal{G}}^{Bess}\left (u\right)
\leq
\pi_{\mathcal{F},\mathcal{G}}\left(u\right),
\end{equation}
for each $u\in E\otimes F$.
\end{enumerate}
\end{proposition} 
\textbf{Proof.}
\begin{enumerate}
\item
Let $u=\sum_{r=1}^{R} x^{r}\otimes y^{r}$ be an arbitrary element of $E\otimes F$. Using the dual norm inequalities \eqref{norms.dual.equi} (from lemma \ref{norm.dual}), the Schauder frame decompositions: 
$$x^{r}=\sum_{m=1}^{\infty}b_{m}^{*}(x^{r})a_{m}\;\text{and}\; y^{r}=\sum_{n=1}^{\infty}y_{n}^{*}(y^{r})x_{n},$$
 and the triangle inequality, we obtain:
\begin{align*}
\mathcal{L}_{\mathcal{F}}^{-1}
\mathcal{L}_{\mathcal{G}}^{-1}
\varepsilon_{\mathcal{F},\mathcal{G}}\left(u\right)
&=
\mathcal{L}_{\mathcal{F}}^{-1}
\mathcal{L}_{\mathcal{G}}^{-1}
\sup_{\substack{\|f^*\|_{E^*,\mathcal{F}}^{Bess} \leq 1 \\ \|g^*\|_{F^*,\mathcal{G}}^{Bess} \leq 1}}
 \left \vert \underset{r=1}{\overset{R}{\sum }}
   f^{*}(x^{r})g^{*}(y^{r})\right\vert\\
     &\leq 
   \mathcal{L}_{\mathcal{F}}^{-1}
\mathcal{L}_{\mathcal{G}}^{-1}
\sup_{\substack{\mathcal{L}_{\mathcal{F}}^{-1}\|f^*\|_{E^*} \leq 1 \\ \mathcal{L}_{\mathcal{G}}^{-1}\|g^*\|_{F^*} \leq 1}}
 \left \vert \underset{r=1}{\overset{R}{\sum}}
   f^{*}(x^{r})g^{*}(y^{r})\right\vert\\
   &\leq  
   Bess\left (\sum_{r=1}^{R}x^{r}\otimes y^{r}\right).
\end{align*}
Taking the infimum over all representations of  $u$ in $E\otimes F$, we obtain
 $$\mathcal{L}_{\mathcal{F}}^{-1}
\mathcal{L}_{\mathcal{G}}^{-1}
\varepsilon_{\mathcal{F},\mathcal{G}}\left (u \right)\leq\alpha_{\mathcal{F},\mathcal{G}}^{Bess}(u).$$
For the right inequality, observe that:
\begin{align*}
\alpha_{\mathcal{F},\mathcal{G}}^{Bess}\left(u\right)
&\leq
\sum_{r=1}^{R}\Vert x^{r}\Vert_{E,\mathcal{F}}^{Bess}\Vert y^{r}\Vert_{E,\mathcal{G}}^{Bess}.
\end{align*}
Taking the infimum over all such representations gives:
$$ \alpha_{\mathcal{F},\mathcal{G}}^{Bess}(u)
\leq
\pi_{\mathcal{F},\mathcal{G}}(u).$$
Hence,
$$\mathcal{L}_{\mathcal{F}}^{-1}
\mathcal{L}_{\mathcal{G}}^{-1}
\varepsilon_{\mathcal{F},\mathcal{G}}\left(u\right)\leq\alpha_{\mathcal{F},\mathcal{G}}^{Bess}(u)
\leq
\pi_{\mathcal{F},\mathcal{G}}(u).$$
\item
\begin{enumerate}
\item
Let $x\otimes y$ be an arbitrary elementary tensor in $E\otimes F$.
By definition of the norms $\Vert\cdot\Vert^{Bess}_{E,\mathcal{F}}$ and $\Vert\cdot\Vert^{Bess}_{F,\mathcal{G}}$, we have:
$$
\Vert x\Vert^{Bess}_{E,\mathcal{F}}\Vert y\Vert^{Bess}_{F,\mathcal{G}}
=
\sup_{\substack{\|f^*\|_{E^*} \leq 1 \\ \|g^*\|_{F^*} \leq 1}}
\underset{m,n=1}{\overset{\infty}{\sum}} 
 \left\vert b_{m}^{*}(x)y_{n}^{*}(y)f^{*}(a_{m})g^{*}\left(x_{n}\right)\right\vert .
 $$\\
Suppose $\sum_{r=1}^{R}x^{r}\otimes y^{r}$ is another representation of 
$x\otimes y$ in $E\otimes F$. Then:
\begin{align*}
\Vert x\Vert^{Bess}_{E,\mathcal{F}}\Vert y\Vert^{Bess}_{F,\mathcal{G}}
&=
\sup_{\substack{\|f^*\|_{E^*} \leq 1 \\ \|g^*\|_{F^*} \leq 1}}
\underset{m,n=1}{\overset{\infty}{\sum}} 
 \left\vert b_{m}^{*}\otimes y_{n}^{*}(x\otimes y)f^{*}(a_{m})g^{*}\left(x_{n}\right)\right \vert\\
&=
\sup_{\substack{\|f^*\|_{E^*} \leq 1 \\ \|g^*\|_{F^*} \leq 1}} 
\underset{m,n=1}{\overset{\infty}{\sum}} 
\left\vert
\underset{r=1}{\overset{R}{\sum}}
b_{m}^{*}(x^{r})y_{n}^{*}(y^{r}) f^{*}(a_{m})g^{*}\left(x_{n}\right)   
\right\vert\\
&\leq
\sup_{\substack{\|f^*\|_{E^*} \leq 1 \\ \|g^*\|_{F^*} \leq 1}}
\underset{r=1}{\overset{R}{\sum}} 
\underset{m,n=1}{\overset{\infty}{\sum}} 
 \left\vert b_{m}^{*}(x^{r})y_{n}^{*}(y^{r})f^{*}(a_{m})g^{*}\left(x_{n}\right)    \right \vert \\
 &=
  Bess\left (\sum_{r=1}^{R}x^{r}\otimes y^{r}\right).  
\end{align*}
Taking the infimum over all such representations of $x\otimes y$ gives:
$$\Vert x\Vert^{Bess}_{E,\mathcal{F}}\Vert y\Vert^{Bess}_{F,\mathcal{G}}\leq\alpha_{\mathcal{F},\mathcal{G}}^{Bess} (x\otimes y).$$
Consequently,
$$\alpha_{\mathcal{F},\mathcal{G}}^{Bess}(x\otimes y)=\Vert x\Vert^{Bess}_{E,\mathcal{F}}\Vert y\Vert^{Bess}_{F,\mathcal{G}}.$$
\item
Let $u=\sum_{r=1}^{R}x^{r}\otimes y^{r}\in E\otimes F$, and let $\lambda\in\mathbb{K}$.  It is clear that: 
$$\alpha_{\mathcal{F},\mathcal{G}}^{Bess}(\lambda u)=\vert\lambda\vert\alpha_{\mathcal{F},\mathcal{G}}^{Bess}(u).$$
\item
Let $u_{1}=\sum_{r_{1}=1}^{R_{1}}x_{1}^{r_{1}}\otimes y_{1}^{r_{1}}$ and
$u_{2}=\sum_{r_{2}=1}^{R_{2}}x_{2}^{r_{2}}\otimes y_{2}^{r_{2}}$ be two elements of $E\otimes F$.  Then:  
\begin{align*}
\alpha_{\mathcal{F},\mathcal{G}}^{Bess}(u_{1}+u_{2})
&\leq
Bess\left (\sum_{r_{1}=1}^{R_{1}}x_{1}^{r_{1}}\otimes y_{1}^{r_{1}}
+\sum_{r_{2}=1}^{R_{2}}x_{2}^{r_{2}}\otimes y_{2}^{r_{2}}\right)\\
&\leq
Bess\left(\sum_{r_{1}=1}^{R_{1}}x_{1}^{r_{1}}\otimes y_{1}^{r_{1}}\right)
+
Bess\left(\sum_{r_{2}=1}^{R_{2}}x_{2}^{r_{2}}\otimes y_{2}^{r_{2}}\right).
\end{align*}
Taking the infimum over all representations of $u_{1}$ and $u_{2}$, yields:
 $$
 \alpha_{\mathcal{F},\mathcal{G}}^{Bess}\left (u_{1}+u_{2}\right)\\
\leq
\alpha_{\mathcal{F},\mathcal{G}}^{Bess}\left (u_{1}\right)
+
\alpha_{\mathcal{F},\mathcal{G}}^{Bess}\left (u_{2}\right).
$$
\item
Let $u\in E\otimes F $ satisfy $\alpha_{\mathcal{F},\mathcal{G}}^{Bess}(u)=0$. From inequality \eqref{New-Norm-Between-inj-proj}, it follows that:
$$u=0.$$
\item  
Let  $f_{0}^{*}\otimes g_{0}^{*}\in E^{*}\otimes F^{*}$ and
$u=\sum_{r=1}^{R}x^{r}\otimes y^{r}\in E\otimes F$. By the Schauder frame decompositions and the triangle inequality, we obtain:
\begin{align*}
\left\vert f_{0}^{*}\otimes g_{0}^{*}\left(u\right)\right\vert
     &\leq 
     \Vert f_{0}^{*}\Vert_{E^{*}}\Vert g_{0}^{*}\Vert _{F^{*}}
Bess\left(\sum_{r=1}^{R}x^{r}\otimes y^{r}\right).
\end{align*}
Taking the infimum over all finite representations of $u$ gives:
\begin{align*}
\left\vert f_{0}^{*}\otimes g_{0}^{*}\left(u\right)\right\vert
     &\leq 
     \Vert f_{0}^{*}\Vert_{E^{*}}\Vert g_{0}^{*}\Vert _{F^{*}}
\alpha_{\mathcal{F},\mathcal{G}}^{Bess}\left(u\right).
\end{align*}
\end{enumerate}
Thus,  $\alpha_{\mathcal{F},\mathcal{G}}^{Bess}$ is a reasonable crossnorm. Consequently, it satisfies:
\begin{equation*}
\varepsilon_{\mathcal{F},\mathcal{G}}\left (u\right)
\leq
\alpha_{\mathcal{F},\mathcal{G}}^{Bess}\left (u\right)
\leq
\pi_{\mathcal{F},\mathcal{G}}\left(u\right),\;(u\in E\otimes F)
\end{equation*}
(see \cite[Theorem 1.1.3, p. 7]{Diestel.2008} or \cite[Proposition 6.1, p. 127]{Raymond.A.Ryan}). This completes the proof.
\end{enumerate}
\qed
\begin{proposition}
\label{norm.for.BSB}
 Suppose that $\mathcal{F}$ and $\mathcal{G}$ are USB for $E$ and $F$, respectively. 
The following holds:
 \begin{enumerate}
 \item
 For each 
$u=\sum_{i,j=1}^{k,k^{'}}\lambda_{i,j}a_{i}\otimes x_{j}$ in the algebraic span of  $\left\lbrace a_{i}\otimes x_{j}: i,j\in \mathbb{N}^{*} \right\rbrace$, we have:
\[
\alpha_{\mathcal{F},\mathcal{G}}^{Bess}\left(u\right)
=
\sup_{\substack{\|f^*\|_{E^*} \leq 1 \\ \|g^*\|_{F^*} \leq 1}}
\underset{i,j=1}{\overset{k,k^{'}}{\sum}} 
 \left\vert\lambda_{i,j}\right \vert
 \left\vert  f^{*}(a_{i})\right \vert 
 \left\vert  g^{*}(x_{j})\right \vert
 .
\]
\item
For each $u$ in $E\otimes F$, we have:
$$u
=
\underset{N\rightarrow\infty}{\lim}
\underset{M\rightarrow\infty}{\lim}
\underset{m,n=1}{\overset{M,N}{\sum}} 
 b_{m}^{*}\otimes y_{n}^{*}\left(u\right)a_{m}\otimes x_{n}$$
in the norm $\alpha_{\mathcal{F},\mathcal{G}}^{Bess}$. Consequently, the linear subspace spanned by 
$$\left\lbrace a_{m}\otimes x_{n}: m,n\in\mathbb{N}^{*}\right\rbrace$$
 is dense in $E\otimes F$.
\end{enumerate}
\end{proposition}
\textbf{Proof.}
\begin{enumerate}
\item
Let $u=\sum_{i,j=1}^{k,k^{'}}\lambda_{i,j}a_{i}\otimes x_{j}$. By the biorthogonality of $\mathcal{F}$ and $\mathcal{G}$:
\begin{align*}
Bess&\left(\sum_{i,j=1}^{k,k^{'}}\lambda_{i,j}a_{i}\otimes x_{j}\right)\\
&=\sup_{\substack{\|f^*\|_{E^*} \leq 1 \\ \|g^*\|_{F^*} \leq 1}}
\underset{i,j=1}{\overset{k,k^{'}}{\sum}} 
\underset{m,n=1}{\overset{\infty}{\sum}} 
 \left\vert
 \lambda_{i,j}b_{m}^{*}(a_{i})y_{n}^{*}(x_{j})f^{*}(a_{m})g^{*}\left(x_{n}\right)    \right \vert \\
 &=
\sup_{\substack{\|f^*\|_{E^*} \leq 1 \\ \|g^*\|_{F^*} \leq 1}}
\underset{i,j=1}{\overset{k,k^{'}}{\sum }} 
 \left\vert \lambda_{i,j}\right \vert
  \left\vert  f^{*}(a_{i})\right \vert 
 \left\vert  g^{*}(x_{j})\right \vert
 .
\end{align*}
Consequently,
$$
\alpha_{\mathcal{F},\mathcal{G}}^{Bess}(u)
\leq
\sup_{\substack{\|f^*\|_{E^*} \leq 1 \\ \|g^*\|_{F^*} \leq 1}}
\underset{i,j=1}{\overset{k,k^{'}}{\sum }} 
 \left\vert\lambda_{i,j}\right\vert
 \left\vert f^{*}(a_{i})\right\vert 
 \left\vert g^{*}(x_{j})\right\vert .
$$
Now, suppose that $\sum_{r=1}^{R}x^{r}\otimes y^{r}$ is another representation of $u$ in $E\otimes F$. Expanding each $x^{r}$ and $y^{r}$ in their respective bases gives:
\begin{align*}
\sum_{i,j=1}^{k,k^{'}}\lambda_{i,j}a_{i}\otimes x_{j}=\sum_{r=1}^{R}x^{r}\otimes y^{r}
&=\sum_{r=1}^{R}
\left (\sum_{i=1}^{\infty}x_{i}^{r}a_{i}\right)\otimes
 \left(\sum_{j=1}^{\infty}y_{j}^{r}x_{j}\right).  
\end{align*}
The biorthogonality relation:
$$b_{m}^{*}\otimes y_{n}^{*}\left(a_{i}\otimes x_{j}\right)=\delta_{m,i}\delta_{n,j},$$ implies that:
 $$\sum_{r=1}^{R}x_i^r y_j^r=\begin{cases} 
\lambda_{i,j} & \text{if } i \leq k \text{ and } j \leq k', \\
0 & \text{otherwise.}
\end{cases}.$$ 
Therefore,
\begin{align*}
Bess\left(\sum_{r=1}^{R}x^{r}\otimes y^{r}\right)
&=
\sup_{\substack{\|f^*\|_{E^*} \leq 1 \\ \|g^*\|_{F^*} \leq 1}}
\underset{r=1}{\overset{R}{\sum}} 
\underset{m,n=1}{\overset{\infty}{\sum}} 
 \left\vert b_{m}^{*}(x^{r})y_{n}^{*}(y^{r})f^{*}(a_{m})g^{*}\left(x_{n}\right)    \right \vert \\
 &=
\sup_{\substack{\|f^*\|_{E^*} \leq 1 \\ \|g^*\|_{F^*} \leq 1}}\underset{r=1}{\overset{R}{\sum}} 
\underset{m,n=1}{\overset{\infty}{\sum}} 
 \left\vert x_{m}^{r}y_{n}^{r} f^{*}(a_{m})g^{*}\left(x_{n}\right)\right\vert\\
  &\geq
\sup_{\substack{\|f^*\|_{E^*} \leq 1 \\ \|g^*\|_{F^*} \leq 1}}\underset{m,n=1}{\overset{\infty}{\sum }} 
 \left\vert
 \underset{r=1}{\overset{R}{\sum }} 
   x_{m}^{r}y_{n}^{r}
    \right \vert
  \left\vert
   f^{*}(a_{m})g^{*}\left(x_{n}\right)\right\vert\\
 &=
\sup_{\substack{\|f^*\|_{E^*} \leq 1 \\ \|g^*\|_{F^*} \leq 1}}\underset{m,n=1}{\overset{k,k^{'}}{\sum }} 
 \left\vert\lambda_{m,n}\right\vert\left\vert f^{*}(a_{m})\right\vert\left\vert g^{*}\left(x_{n}\right)\right\vert 
  .
\end{align*}
Taking the infimum over all representations yields the reverse inequality, hence the equality.
\item
Let $x\otimes y\in E\otimes F$. Define the finite-rank projections:
$$
P_{M}(x)=\sum_{m=1}^{M}b_{m}^{*}(x)a_{m},\;\;\;Q_{N}(y)=\sum_{n=1}^{N}y_{n}^{*}(y)x_{n}.
$$
The triangle inequality yields:
\begin{align*}
\alpha_{\mathcal{F},\mathcal{G}}^{Bess}&\left(x\otimes y-\sum_{m,n=1}^{M,N}b_{m}^{*}(x)y_{n}^{*}(y)a_{m}\otimes x_{n}\right)\\
&=
\alpha_{\mathcal{F},\mathcal{G}}^{Bess}\left(x\otimes y-P_{M}(x)\otimes Q_{N}(y)\right)\\
&=
\alpha_{\mathcal{F},\mathcal{G}}^{Bess}\left((x-P_{M}(x))\otimes y+P_{M}(x)\otimes (y-Q_{N}(y))\right)\\
&\leq
\alpha_{\mathcal{F},\mathcal{G}}^{Bess}\left((x-P_{M}(x))\otimes y\right)
+\alpha_{\mathcal{F},\mathcal{G}}^{Bess}\left(P_{M}(x)\otimes (y-Q_{N}(y))\right)\\
&=
\Vert x-P_{M}(x)\Vert_{E,\mathcal{F}}^{Bess}\Vert y\Vert^{Bess}_{F,\mathcal{G}}
+
\Vert P_{M}(x)\Vert_{E,\mathcal{F}}^{Bess}\Vert y-Q_{N}(y)\Vert^{Bess}_{F,\mathcal{G}}.
\end{align*}
Since $\mathcal{F}$ and $\mathcal{G}$ are SB,  the right-hand side tends to $0$ as
$M,N\rightarrow \infty$. Thus, 
 $$x\otimes y
\stackrel{\alpha_{\mathcal{F},\mathcal{G}}^{Bess}}{=} 
\underset{N\rightarrow\infty}{\lim}
\underset{M\rightarrow\infty}{\lim}
\underset{m,n=1}{\overset{M,N}{\sum}} 
 b_{m}^{*}(x)y_{n}^{*}(y)a_{m}\otimes x_{n}.$$
By linearity, this extends to finite sums $u=\sum_{r=1}^{R}x^{r}\otimes y^{r}$:
$$u
\stackrel{\alpha_{\mathcal{F},\mathcal{G}}^{Bess}}{=} 
\underset{N\rightarrow\infty}{\lim}
\underset{M\rightarrow\infty}{\lim}
\underset{m,n=1}{\overset{M,N}{\sum}} 
 b_{m}^{*}\otimes y_{n}^{*}\left(u\right)a_{m}\otimes x_{n}.$$
This implies  the density of $\text{span}\left\lbrace a_{i}\otimes x_{j}: i,j\in \mathbb{N}^{*} \right\rbrace$ in $\left(E\otimes F,\alpha_{\mathcal{F},\mathcal{G}}^{Bess}\right)$.  
Which finishes the proof.
\end{enumerate}
\qed\\

In the following, we establish that the Besselian crossnorm $\alpha_{\mathcal{F},\mathcal{G}}^{Bess}$  satisfies the second condition of the unconditionality test. Furthermore, we present illustrative examples of tensor product spaces to facilitate a deeper understanding of computations involving this new crossnorm.
\begin{proposition}
\label{examples-classical-spaces}
Let $\mathcal{F}=\mathcal{G}=\left(\left(e_{m},e_{m}^{*}\right)\right)_{m\in \mathbb{N}^{*}}$ be the canonical USB for the classical Banach spaces $\ell_{p}\;(1\leq p<\infty)$ or $c_{0}$. Then the Besselian crossnorm $\alpha_{\mathcal{F},\mathcal{G}}^{Bess}$ has the following properties:
\begin{enumerate} 
\item
On $c_{0}\otimes c_{0}$, $\alpha_{\mathcal{F},\mathcal{G}}^{Bess}=\varepsilon_{\mathcal{F},\mathcal{G}}=\varepsilon$. Consequently,
$$
c_{0}\otimes_{\alpha_{\mathcal{F},\mathcal{G}}^{Bess}} c_{0}
=c_{0}\otimes_{\varepsilon} c_{0}
\simeq c_{0}
.$$ 
\item
On $\ell_{1}\otimes \ell_{1}$, $\alpha_{\mathcal{F},\mathcal{G}}^{Bess}=\pi_{\mathcal{F},\mathcal{G}}=\pi$. Consequently,
$$
\ell_{1}\otimes_{\alpha_{\mathcal{F},\mathcal{G}}^{Bess}} \ell_{1}
=\ell_{1}\otimes_{\pi} \ell_{1}
\simeq \ell_{1} 
.$$ 
\item
On $\ell_{2}\otimes \ell_{2}$, $\alpha_{\mathcal{F},\mathcal{G}}^{Bess}$ is not proportional to either $\varepsilon_{\mathcal{F},\mathcal{G}}$ or $\pi_{\mathcal{F},\mathcal{G}}$.
\item
$\ell_{p}\otimes_{\alpha_{\mathcal{F},\mathcal{G}}^{Bess}} \ell_{1}\simeq\ell_{p}(\ell_{1})$.
\item
$c_{0}\otimes_{\alpha_{\mathcal{F},\mathcal{G}}^{Bess}}F\simeq c_{0}(F)$, if $\mathcal{G}=\left(\left(x_{n},y_{n}^{*}\right)\right)_{n\in \mathbb{N}^{*}}$ is a USB for $F$.
\end{enumerate}
\end{proposition}
\textbf{Proof} 
\begin{enumerate}
\item
Let $u=\sum_{i,j=1}^{k,k^{'}}\lambda_{i,j} e_{i}\otimes e_{j}\in c_{0}\otimes c_{0}$. From example \eqref{example.norm.bess=norm}, we have:
$$\Vert f^{*}\Vert_{c_{0}^{*},\mathcal{F}}^{Bess}=\Vert f^{*}\Vert_{c_{0}^{*}}$$
for each $f^{*}\in c_{0}^{*}\simeq\ell_{1}$. Then,
\begin{align*}
\varepsilon_{\mathcal{F},\mathcal{G}}\left(u\right)
 &:=
 \sup_{\substack{\|f^*\|_{\ell_{1},\mathcal{F}}^{Bess} \leq 1 \\ \|g^*\|_{\ell_{1},\mathcal{G}}^{Bess} \leq 1}}
 \left\vert
\underset{i,j=1}{\overset{k,k^{'}}{\sum }} 
 \lambda_{i,j}f^{*}(e_{i})g^{*}(e_{j})\right \vert\\
 &=
 \sup_{\substack{\|f^*\|_{\ell_{1}} \leq 1 \\ \|g^*\|_{\ell_{1}} \leq 1}}
 \left\vert
\underset{i,j=1}{\overset{k,k^{'}}{\sum }} 
 \lambda_{i,j}f^{*}(e_{i})g^{*}(e_{j})\right \vert\\
&\leq
 \underset{i,j}{\max} \vert \lambda_{i,j}\vert
\sup_{\substack{\|f^*\|_{\ell_{1}} \leq 1 \\ \|g^*\|_{\ell_{1}} \leq 1}}
\underset{i,j=1}{\overset{k,k^{'}}{\sum }} 
 \left\vert f^{*}(e_{i})\right \vert \left\vert g^{*}(e_{j})\right\vert\\
 &=
 \underset{i,j}{\max} \vert \lambda_{i,j}\vert
 \sup_{\substack{\|f^*\|_{\ell_{1}} \leq 1 \\ \|g^*\|_{\ell_{1}} \leq 1}}
\left (
\underset{i=1}{\overset{k}{\sum}} 
 \left\vert f^{*}(e_{i})\right\vert 
 \right )
 \left (
 \underset{j=1}{\overset{k^{'}}{\sum}} 
 \left\vert g^{*}(e_{j})\right\vert
\right)\\
&\leq
 \underset{i,j}{\max}\vert\lambda_{i,j}\vert .
\end{align*}
To show equality, suppose $\underset{i,j}{\max} \vert \lambda_{i,j}\vert=\vert\lambda_{i_{0},j_{0}}\vert$. Choose 
$$f^{*}=e_{i_{0}}^{*} \;\text{and}\;g^{*}=e_{j_{0}}^{*},$$ 
(the coordinate functionals). Then
$$
\left\vert
\underset{i,j=1}{\overset{k}{\sum}} 
 \lambda_{i,j}f^{*}(e_{i})g^{*}(e_{j})\right\vert 
 =\\
 \underset{i,j}{\max}\vert\lambda_{i,j}\vert .
 $$
Hence,
 $$
 \varepsilon_{\mathcal{F},\mathcal{G}}\left(u\right)
=
\underset{i,j}{\max}\vert\lambda_{i,j}\vert .
$$
On the other hand, by proposition \eqref{norm.for.BSB}, it follows that:
\begin{align*}
\alpha_{\mathcal{F},\mathcal{G}}^{Bess}\left(u\right)
&=
\sup_{\substack{\|f^*\|_{\ell_{1}} \leq 1 \\ \|g^*\|_{\ell_{1}} \leq 1}}
\underset{i,j=1}{\overset{k,k^{'}}{\sum }} 
 \left\vert
 \lambda_{i,j}f^{*}(e_{i})g^{*}(e_{j})\right\vert\\
 &\leq
 \underset{i,j}{\max}\vert\lambda_{i,j}\vert
\sup_{\substack{\|f^*\|_{\ell_{1}} \leq 1 \\ \|g^*\|_{\ell_{1}} \leq 1}}\left (
\underset{i=1}{\overset{k}{\sum }} 
 \left\vert  f^{*}(e_{i})\right\vert 
 \right)
 \left(
 \underset{j=1}{\overset{k^{'}}{\sum }} 
 \left\vert g^{*}(e_{j})\right \vert
\right )\\
&\leq
\underset{i,j}{\max} \vert \lambda_{i,j}\vert 
=
\varepsilon_{\mathcal{F},\mathcal{G}}(u).
\end{align*}
Equality \eqref{2-New-Norm-Between-inj-proj} then implies:
 $$
 \alpha_{\mathcal{F},\mathcal{G}}^{Bess}\left(u\right)
=
 \varepsilon_{\mathcal{F},\mathcal{G}}\left (u \right).
$$
By the known isometric isomorphism $c_{0}\otimes_{\varepsilon}c_{0}\simeq c_{0}$ (See \cite[p. 68]{Raymond.A.Ryan}), it follows that:
$$c_{0}\otimes_{\alpha_{\mathcal{F},\mathcal{G}}^{Bess}}c_{0}\simeq c_{0}.$$
This completes the proof of Part 1.  
\item
Let $u=\sum_{i,j=1}^{k,k^{'}} \lambda_{i,j} e_{i}\otimes e_{j}\in \ell_{1}\otimes \ell_{1}$. By proposition \eqref{norm.for.BSB}, we have:
\begin{align*}
\alpha_{\mathcal{F},\mathcal{G}}^{Bess}\left(u\right)
&=
\sup_{\substack{\|f^*\|_{\ell_{\infty}} \leq 1 \\ \|g^*\|_{\ell_{\infty}} \leq 1}}
\underset{i,j=1}{\overset{k,k^{'}}{\sum}} 
 \left\vert\lambda_{i,j}\right\vert\left\vert f^{*}(e_{i})\right\vert\left\vert g^{*}(e_{j})\right \vert
 =
 \underset{i,j=1}{\overset{k,k^{'}}{\sum}} 
 \left\vert\lambda_{i,j}\right\vert
 =\pi(u) .
\end{align*}
By the known isometric isomorphism $\ell_{1}\otimes_{\pi}\ell_{1}\simeq \ell_{1}$ (See, \cite[p. 43]{Raymond.A.Ryan}), it follows that:
$$\ell_{1}\otimes_{\alpha_{\mathcal{F},\mathcal{G}}^{Bess}}\ell_{1}\simeq \ell_{1}.$$
This completes the proof of Part 2. 
\item
Note that for any elementary tensor $x\otimes y\in \ell_{2}\otimes \ell_{2}$, we have:
$$
\varepsilon_{\mathcal{F},\mathcal{G}}(x\otimes y)
=
\alpha_{\mathcal{F},\mathcal{G}}^{Bess}(x\otimes y)
=
\pi_{\mathcal{F},\mathcal{G}}(x\otimes y)
=
\Vert x\Vert_{\ell_{2}}\Vert y\Vert_{\ell_{2}}.
$$
To show that these norms are not proportional, it suffices to demonstrate the existence of:
\begin{itemize}
\item
An element $u\in\ell_{2}\otimes\ell_{2}$ satisfying 
$\alpha_{\mathcal{F},\mathcal{G}}^{Bess}(u) \neq \varepsilon_{\mathcal{F},\mathcal{G}}(u)$,
\item
An element $v\in\ell_{2}\otimes\ell_{2}$ satisfying 
$\alpha_{\mathcal{F},\mathcal{G}}^{Bess}(v) \neq \pi_{\mathcal{F},\mathcal{G}}(v)$.
\end{itemize}
\begin{enumerate}
\item
Let $u=e_{1}\otimes e_{1}+e_{2}\otimes e_{2}$ be an element of $\ell_{2}\otimes \ell_{2}$. Then,
\begin{align*}
\alpha_{\mathcal{F},\mathcal{G}}^{Bess}(u)
&=
\sup_{\substack{\|f^*\|_{\ell_{2}} \leq 1 \\ \|g^*\|_{\ell_{2}} \leq 1}}
\left(
\left\vert  f^{*}(e_{1})\right \vert \left\vert g^{*}(e_{1}) \right \vert
+
 \left\vert f^{*}(e_{2})\right \vert \left\vert g^{*}(e_{2}) \right \vert\right).
\end{align*}
Since the supremum depends only on the first two coordinates, we may restrict to functionals of the form:
$$
f^{*}=(f_{1}^{*},f_{2}^{*},0,0,...)\;\;\text{and}\;\;g^{*}=(g_{1}^{*},g_{2}^{*},0,0...),
$$
where 
$$\vert f_{1}^{*}\vert^{2}+\vert f_{2}^{*}\vert^{2}\leq 1\;\;\text{and}
\;\;\vert g_{1}^{*}\vert^{2}+\vert g_{2}^{*}\vert^{2}\leq 1.$$
Using polar coordinates, write:
$$
f^{*}=(sin(\alpha),cos(\alpha),0,0,...)\;\;\text{and}\;\;g^{*}=(sin(\beta),cos(\beta),0,0...),
$$
where $0\leq \alpha ,\beta\leq\pi/2$. We obtain:
\begin{align*}
\alpha_{\mathcal{F},\mathcal{G}}^{Bess}(u)
&=
\underset{0\leq \alpha ,\beta\leq\pi/2}{\sup}
\left(sin(\alpha)sin(\beta)
 +
 cos(\alpha)cos(\beta)\right )\\
 &=
 \underset{0\leq\alpha ,\beta\leq\pi/2}{\sup}
cos(\alpha-\beta)=1.
\end{align*}
On the other hand,
$\pi_{\mathcal{F},\mathcal{G}}(u)=2$ (See \cite[p. 23]{Raymond.A.Ryan}).
This shows that $\alpha_{\mathcal{F},\mathcal{G}}^{Bess}$ and $\pi_{\mathcal{F},\mathcal{G}}$ are not proportional.
\item
Let $v=e_{1}\otimes e_{1}+e_{1}\otimes e_{2}
+e_{2}\otimes e_{2}-e_{2}\otimes e_{1}\in \ell_{2}\otimes \ell_{2}$. Define:
$$
A_{f^{*},g^{*}}^{v}:=f^{*}(e_{1})g^{*}(e_{1})
   + 
   f^{*}(e_{1})g^{*}(e_{2})
   +
   f^{*}(e_{2})g^{*}(e_{2})
   -
   f^{*}(e_{2})g^{*}(e_{1}),
$$
$$
B_{\alpha,\beta}:=sin(\alpha) sin(\beta)
   + 
   sin(\alpha)cos(\beta)
   +
   cos(\alpha)cos(\beta)
   -
   cos(\alpha)sin(\beta).
$$
Then,
\begin{align*}
\varepsilon_{\mathcal{F},\mathcal{G}}\left (v \right)
&=
\sup_{\substack{\|f^*\|_{\ell_{2}} \leq 1 \\ \|g^*\|_{\ell_{2}} \leq 1}}\left\vert
A_{f^{*},g^{*}}^{v}
   \right \vert
=\underset{-\pi<\alpha , \beta\leq\pi}{\sup}
  \left\vert B_{\alpha,\beta}
    \right \vert \\
    &=   \underset{-\pi<\alpha , \beta\leq\pi}{\sup}
  \left\vert
   cos(\alpha-\beta)
   + 
   sin(\alpha-\beta)
    \right \vert \\
    &=   \underset{-\pi<\alpha , \beta\leq\pi}{\sup}
  \sqrt{2}
  \left\vert
   sin(\alpha-\beta+\pi/4)
    \right \vert
    =\sqrt{2}.
\end{align*}
On the other hand,
\begin{align*}
\alpha_{\mathcal{F},\mathcal{G}}^{Bess}\left (v \right)
   &=
\sup_{\substack{\|f^*\|_{\ell_{2}} \leq 1 \\ \|g^*\|_{\ell_{2}} \leq 1}}\left(\vert f^{*}(e_{1})\vert+\vert f^{*}(e_{2})\vert\right)
   \left(\vert g^{*}(e_{1})\vert +\vert g^{*}(e_{2})\vert\right) 
   \\
&=   \underset{-\pi<\alpha , \beta\leq\pi}{\sup}
(sin(\alpha)+ cos(\alpha))
(sin(\beta)+cos(\beta)) \\
    &=   \underset{-\pi<\alpha , \beta\leq\pi}{\sup}
   2sin(\alpha +\pi/4)
    sin( \beta+\pi/4)
    =2.
\end{align*}
Then $\alpha_{\mathcal{F},\mathcal{G}}^{Bess}$ and $\varepsilon_{\mathcal{F},\mathcal{G}}$ are not proportional.
\end{enumerate}
This complete the proof of Part 3.
\item 
We first prove that the linear operator
\begin{align*}
J:\;\;\;\;  \ell_{p}\otimes \ell_{1}&\longrightarrow  \ell_{p}(\ell_{1}) \\ 
 \sum_{i,j=1}^{k,k^{'}} \lambda_{i,j} e_{i}\otimes e_{j}&  \longmapsto  
 \left(\sum_{j=1}^{k'}\lambda_{i,j}e_{j}\right )_{1\leq i \leq k}
\end{align*}
is an isometric isomorphism. \\
Let $u=\sum_{i,j=1}^{k,k^{'}} \lambda_{i,j} e_{i}\otimes e_{j}\in\ell_{p}\otimes \ell_{1}$. Then, 
\begin{align*}
\alpha_{\mathcal{F},\mathcal{G}}^{Bess}\left (u \right)
  &=
  \sup_{\substack{\|f^*\|_{\ell_{p^{*}}} \leq 1 \\ \|g^*\|_{\ell_{\infty}} \leq 1}}
\underset{i=1}{\overset{k}{\sum}}
\left (  
\underset{j=1}{\overset{k^{'}}{\sum}}
 \left\vert \lambda_{i,j}\right\vert
 \left\vert g^{*}(e_{j})\right\vert 
 \right )
 \left\vert f^{*}(e_{i})\right\vert\\
 &=
\underset{\left\Vert f^{*}\right\Vert_{\ell_{p^{*}}}\leq 1}{\sup}
\underset{i=1}{\overset{k}{\sum }}
\left ( 
\underset{j=1}{\overset{k^{'}}{\sum}} 
 \left\vert\lambda_{i,j}\right\vert
 \right)
 \left\vert f^{*}(e_{i})\right\vert
 = 
 \left(
\underset{i=1}{\overset{k}{\sum}}
 \left\Vert\sum_{j=1}^{j=k'}\lambda_{i,j}e_{j}\right\Vert_{\ell_{1}}^{p}
 \right )^{1/p}\\
 &=\Vert J(u)\Vert_{\ell_{p}(\ell_{1})}
 .
\end{align*}
Then $J$ is an isometric isomorphism. By the Hahn-Banach theorem, $J$
extends uniquely to an isometric operator
$$J:\ell_{p}\otimes_{\alpha_{\mathcal{F},\mathcal{G}}^{Bess}}\ell_{1}\rightarrow \ell_{p}(\ell_{1}).$$
To show  surjectivity, let $(y^{i})_{i\in\mathbb{N}^{*}}\in \ell_{p}(\ell_{1})$. 
For any integers $M\leq N$,
\begin{align*}
\alpha_{\mathcal{F},\mathcal{G}}^{Bess}\left(\sum_{ i=M}^{N}e_{i}\otimes y^{i}\right)
&\leq
\sup_{\substack{\|f^*\|_{\ell_{p^{*}}} \leq 1 \\ \|g^*\|_{\ell_{\infty}} \leq 1}}
\sum_{i=M}^{N}
\underset{m,n=1}{\overset{\infty}{\sum }} 
 \left\vert e_{m}^{*}(e_{i})e_{n}^{*}(y^{i})f^{*}(e_{m})g^{*}\left(e_{n}\right)    \right \vert \\
 &=
\sup_{\substack{\|f^*\|_{\ell_{p^{*}}} \leq 1 \\ \|g^*\|_{\ell_{\infty}} \leq 1}}
\sum_{i=M}^{N}
\underset{n=1}{\overset{\infty}{\sum}} 
 \left\vert f^{*}(e_{i})\right\vert
  \left\vert e_{n}^{*}(y^{i})g^{*}\left(e_{n}\right)\right\vert \\
  &\leq
\underset{\left \Vert f^{*}\right\Vert_{\ell_{p^{*}}}\leq 1}{\sup} 
\sum_{i=M}^{N}
\left\vert   f^{*}(e_{i})\right \vert
\underset{\left\Vert g^{*}\right\Vert_{\ell_{\infty}}\leq 1}{\sup} 
\underset{n=1}{\overset{\infty}{\sum }} 
  \left\vert e_{n}^{*}(y^{i})g^{*}\left(e_{n}\right)\right\vert\\
  &\leq
\underset{\left \Vert f^{*}\right\Vert_{\ell_{p^{*}}} \leq 1}{\sup} 
\sum_{i=M}^{N}
\left\vert   f^{*}(e_{i})\right \vert
  \left\Vert y^{i}\right\Vert_{\ell_{1}}
=
\left (
\sum_{i=M}^{N}
  \left\Vert y_{i}\right\Vert_{\ell_{1}}^{p}\right )^{1/p}  .
\end{align*}
Therefore, $\sum_{i}e_{i}\otimes y^{i}$ is a Cauchy series in the Banach space
$\ell_{p}\otimes_{\alpha_{\mathcal{F},\mathcal{G}}^{Bess}}\ell_{1}$,  thus converges to an element in
$u\in\ell_{p}\otimes_{\alpha_{\mathcal{F},\mathcal{G}}^{Bess}}\ell_{1}$.  Since 
$$
J\left (\sum_{ i=1}^{\infty}e_{i}\otimes y^{i}\right)=(y^{i})_{i\in\mathbb{N}^{*}},$$
then  $J$ is surjective. This completes the proof of Part 4.
\item
Define the linear map:
\begin{align*}
J:\;\;\;\;  c_0 \otimes F &\longrightarrow  c_0(F)\\ 
 \sum_{i,j=1}^{k,k^{'}}\lambda_{i,j} e_{i}\otimes x_{j}&  \longmapsto  
 \left(\sum_{j=1}^{k'}\lambda_{i,j}x_{j}\right )_{1\leq i \leq k}
\end{align*}
Our goal is to show that \( J \) is an isometry when \( c_0 \otimes E \) is equipped with the norm \( \alpha_{\mathcal{F},\mathcal{G}}^{\text{Bess}} \), and that \( J \) extends to an isometric isomorphism on the completions.\\
Let \( u=\sum_{i,j=1}^{k,k^{'}}\lambda_{i,j} e_{i}\otimes x_{j}\in c_0 \otimes F\), we have:
\begin{align*}
\alpha_{\mathcal{F},\mathcal{G}}^{Bess}\left (u \right)
  &=
  \sup_{\substack{\|f^*\|_{\ell_{1}} \leq 1 \\ \|g^*\|_{F^{*}} \leq 1}}
\underset{i=1}{\overset{k}{\sum}}
\left (  
\underset{j=1}{\overset{k^{'}}{\sum}}
 \left\vert \lambda_{i,j}\right\vert
 \left\vert g^{*}(x_{j})\right\vert 
 \right )
 \left\vert f^{*}(e_{i})\right\vert\\
 &=
\underset{\left\Vert g^{*}\right\Vert_{F^{*}}\leq 1}{\sup}
\underset{1\leq i\leq k}{\max}
\underset{j=1}{\overset{k^{'}}{\sum}}
 \left\vert \lambda_{i,j}\right\vert
 \left\vert g^{*}(x_{j})\right\vert \\
&=
\underset{1\leq i\leq k}{\max}
\underset{\left\Vert g^{*}\right\Vert_{F^{*}}\leq 1}{\sup}
\underset{j=1}{\overset{k^{'}}{\sum}}
 \left\vert \lambda_{i,j}\right\vert
 \left\vert g^{*}(x_{j})\right\vert
=
\underset{1\leq i\leq k}{\max}
\left\Vert
\underset{j=1}{\overset{k^{'}}{\sum}}
 \lambda_{i,j}x_{j}\right\Vert_{F}\\
 &=\left\Vert J(u)\right\Vert_{c_{0}(F)}
\end{align*}
Then $J$ is an isometric isomorphism. By the Hahn-Banach theorem, $J$
extends uniquely to an isometric operator
$$J:c_{0}\otimes_{\alpha_{\mathcal{F},\mathcal{G}}^{Bess}}F\rightarrow c_{0}(F).$$
To show  surjectivity, let $(y^{i})_{i\in\mathbb{N}^{*}}\in c_{0}(F)$. 
For any integers $M\leq N$,
\begin{align*}
\alpha_{\mathcal{F},\mathcal{G}}^{Bess}\left(\sum_{ i=M}^{N}e_{i}\otimes y^{i}\right)
&\leq
\sup_{\substack{\|f^*\|_{\ell_{1}} \leq 1 \\ \|g^*\|_{F^{*}} \leq 1}}
\sum_{i=M}^{N}
\underset{m,n=1}{\overset{\infty}{\sum }} 
 \left\vert e_{m}^{*}(e_{i})y_{n}^{*}(y^{i})f^{*}(e_{m})g^{*}\left(x_{n}\right)    \right \vert \\
 &=
\sup_{\substack{\|f^*\|_{\ell_{1}} \leq 1 \\ \|g^*\|_{F^{*}} \leq 1}}
\sum_{i=M}^{N}
\underset{n=1}{\overset{\infty}{\sum}} 
 \left\vert f^{*}(e_{i})\right\vert
  \left\vert y_{n}^{*}(y^{i})g^{*}\left(x_{n}\right)\right\vert \\
  &\leq
\underset{\left \Vert f^{*}\right\Vert_{\ell_{1}}\leq 1}{\sup} 
\sum_{i=M}^{N}
\left\vert   f^{*}(e_{i})\right \vert
\underset{\left\Vert g^{*}\right\Vert_{F^{*}}\leq 1}{\sup} 
\underset{n=1}{\overset{\infty}{\sum }} 
  \left\vert y_{n}^{*}(y^{i})g^{*}\left(x_{n}\right)\right\vert\\
  &=
\underset{\left \Vert f^{*}\right\Vert_{\ell_{1}} \leq 1}{\sup} 
\sum_{i=M}^{N}
\left\vert   f^{*}(e_{i})\right \vert
  \left\Vert y^{i}\right\Vert_{F,\mathcal{G}}^{Bess}
=
\underset{M\leq i\leq N}{\max}
\left\Vert y^{i}\right\Vert_{F,\mathcal{G}}^{Bess}
\end{align*}
Therefore, $\sum_{i}e_{i}\otimes y^{i}$ is a Cauchy series in the Banach space
$c_{0}\otimes_{\alpha_{\mathcal{F},\mathcal{G}}^{Bess}}F$,  thus converges to an element in
$u\in c_{0}\otimes_{\alpha_{\mathcal{F},\mathcal{G}}^{Bess}}F$.  Since 
$$
J\left (\sum_{ i=1}^{\infty}e_{i}\otimes y^{i}\right)=(y^{i})_{i\in\mathbb{N}^{*}},$$
then  $J$ is surjective. This completes the proof Part 5.
\end{enumerate}
\qed
\begin{remark}
In general, the Besselian crossnorm $\alpha_{\mathcal{F},\mathcal{G}}^{Bess}$ is not uniform, despite coinciding with the injective tensor norm on $c_{0}\otimes c_{0}$ and with the projective tensor norm on $\ell_{1}\otimes \ell_{1}$ (both of which are uniform norms).
\end{remark}
\textbf{Proof.}
Consider $E=F=\ell_{2}$ with the canonical unconditional Schauder basis $\mathcal{F}=\mathcal{G}=\left(\left(e_{m},e_{m}^{*}\right)\right)_{m\in \mathbb{N}^{*}}$.  The Besselian constants satisfy $\mathcal{L}_{\mathcal{F}}=\mathcal{L}_{\mathcal{G}}=1$. Define the following:
\begin{itemize}
\item
A bounded linear operator $S \in L(\ell_{2})$ by: 
$$S(e_1)=S(e_2)=e_1,\,S(e_{i})=0,\;\text{for}\; i\geq 3.$$
The norm of $S$ is $\|S\|_{L(\ell_{2})} = \sqrt{2}$.
\item 
Let $T = I$ be the identity operator on $\ell_{2}$, with $\|T\|_{L(\ell_{2})} = 1$.
\item 
Define the tensor $u = e_1 \otimes e_1 + e_2 \otimes e_2\in \ell_{2}\otimes\ell_{2}$. The proof of proposition \eqref{examples-classical-spaces} implies that:
\end{itemize}
\[
\alpha_{\mathcal{F},\mathcal{G}}^{\text{Bess}}(u)=1
\]
Now consider the tensor $v=S\otimes T(u)$:
\[
v = S(e_1) \otimes T(e_1) + S(e_2) \otimes T(e_2) = e_1 \otimes e_1 + e_1 \otimes e_1 = 2e_1 \otimes e_1
\]
The Besselian crossnorm of $v$ is:
\[
\alpha_{\mathcal{F},\mathcal{G}}^{\text{Bess}}(v) = \sup_{\substack{\|f^*\|_{\ell_{2}} \leq 1 \\ \|g^*\|_{\ell_{2}} \leq 1}} 2|f^*(e_1)| |g^*(e_1)| = 2.
\]
Then
\[
\alpha_{\mathcal{F},\mathcal{G}}^{\text{Bess}}(v) = 2 > \sqrt{2}= \|S\|_{L(\ell_{2})}\|T\|_{L(\ell_{2})} \alpha_{\mathcal{F},\mathcal{G}}^{\text{Bess}}(u).
\]
Violation of the uniformity condition. Hence, $\alpha_{\mathcal{F},\mathcal{G}}^{\text{Bess}}$ is not uniform.\qed
\begin{remark}
\label{norm.signs.2}
Suppose that $\mathcal{F}$ and $\mathcal{G}$ are USB for 
$\left(E,\Vert\cdot\Vert_{E}\right)$ and $\left (F,\Vert\cdot\Vert_{F}\right)$, respectively. We have,
\begin{align*}
\alpha_{\mathcal{F},\mathcal{G}}^{Bess}\left(\sum_{i,j=1}^{k,k'}\lambda_{i,j}a_{i}\otimes x_{j}\right)
&=
\alpha_{\mathcal{F},\mathcal{G}}^{Bess}\left (\sum_{i,j=1}^{k,k'}\varepsilon_{i,j}a_{i}\otimes x_{j}\right)\\
&=
\alpha_{\mathcal{F},\mathcal{G}}^{Bess}\left (\sum_{i,j=1}^{k,k'}\vert\lambda_{i,j}\vert a_{i}\otimes x_{j}\right), 
\end{align*}
for each element $\sum_{i,j=1}^{k,k'}\lambda_{i,j}a_{i}\otimes x_{j}$ of $span\left\lbrace a_{i}\otimes x_{j} \right\rbrace$, and each sequence of signs $(\varepsilon_{m})_{m\in\mathbb{N}^{*}}\in \mathcal{S}$.
\end{remark}
\begin{theorem}\cite[Theorem 18.1, p. 172]{I.Singer.I}
Suppose $\mathcal{F}$ and $\mathcal{G}$ are Schauder bases for $E$ and $F$, respectively, and suppose $\alpha$ is a uniform reasonable crossnorm on $E\otimes F$.
Then, the tensor product sequence
$$\mathcal{F}\otimes\mathcal{G}=\left((a_{m}\otimes x_{n},b_{m}^{*}\otimes y_{n}^{*})\right)_{m,n\in\mathbb{N}^{*}},$$
ordered in the square ordering, forms a Schauder bases for $E\otimes_{\alpha}F$.
\end{theorem}

In this theorem, the uniformity of $\alpha$ is essential to guarantee that $\mathcal{F}\otimes\mathcal{G}$ is a Schauder basis for $E\otimes_{\alpha}F$. In contrast, in our setting, we prove that $\mathcal{F}\otimes\mathcal{G}$ is a USB for  $E\otimes_{\alpha_{\mathcal{F},\mathcal{G}}^{Bess}}F$ although $\alpha_{\mathcal{F},\mathcal{G}}^{Bess}$ is  generally not uniform. To achieve this goal, we use the following characterization:
\begin{theorem} \cite[Theorem 16.1, p. 461]{I.Singer.I}.
\label{characterization.USB.Singer}
Let $G$ be a Banach space and $\left((z_{k},h_{k}^{*})\right)_{k\in\mathbb{N}^{*}}$ be a sequence in $G\times G^{*}$ such that:
the set of all finite linear combinations 
$$\left\lbrace\sum_{k=1}^{n}\lambda_{k}z_{k},\,\,\,\lambda_{k}\in \mathbb{K}, k=1,...,n; n=1,2,...\right\rbrace$$
 is dense in $G$, and the sequence $\left((z_{k},h_{k}^{*})\right)_{k\in\mathbb{N}^{*}}$ is a biorthogonal system.
Then the following are equivalent:
\begin{enumerate}
\item
The sequence $\left((z_{k},h_{k}^{*})\right)_{k\in\mathbb{N}^{*}}$ is an unconditional Schauder basis for $G$.
\item
For each $z\in G$ and $h^{*}\in G^{*}$ the numerical series
$$
\sum_{k=1}^{\infty}\vert h_{k}^{*}(z)h^{*}(z_{k})\vert
$$
converges.
\end{enumerate}
\end{theorem}
\begin{theorem}
\label{tensor.bess.is.bess}
Suppose $\mathcal{F}$ and $\mathcal{G}$ are USB for $E$ and $F$ respectively.
Let $\alpha_{\mathcal{F},\mathcal{G}}^{Bess}$ be the Besselian crossnorm associated with $\mathcal{F}$ and $\mathcal{G}$.
Then, the tensor product sequence
$$\mathcal{F}\otimes\mathcal{G}=\left((a_{m}\otimes x_{n},b_{m}^{*}\otimes y_{n}^{*})\right)_{m,n\in\mathbb{N}^{*}},$$
ordered in the square ordering, forms a USB for $E\otimes_{\alpha_{\mathcal{F},\mathcal{G}}^{Bess}}F$.
\end{theorem}
\textbf{Proof.}
Theorem \eqref{unc.impl.bess} implies that the sequences $\mathcal{F}$ and
$\mathcal{G}$ are BSB for $E$ and $F$, respectively. Then, by proposition \ref{New-crossnorm}, $\alpha_{\mathcal{F},\mathcal{G}}^{Bess}$ is a reasonable crossnorm.
Let  $A^{*}\in\left(E\otimes_{\alpha_{\mathcal{F},\mathcal{G}}^{Bess}}F\right )^{*}$ and
$u\in E\otimes_{\alpha_{\mathcal{F},\mathcal{G}}^{Bess}}F$.  Take a sequence $(u_{k})_{k\in\mathbb{N}^{*}}$ in $E\otimes F$ such that:
$$\lim_{k\rightarrow \infty}u_{k}=u.$$ 
Let $k,M,N\in\mathbb{N}^{*}$. For each 
$1\leq m\leq M$ and $1\leq n\leq N$, define the sign factor $\varepsilon_{m,n}$ such that:  
$$\vert b_{m}^{*}\otimes y_{n}^{*}(u_{k})A^{*}(a_{m}\otimes x_{n})\vert
=
\varepsilon_{m,n}b_{m}^{*}\otimes y_{n}^{*}(u_{k})A^{*}(a_{m}\otimes x_{n}).$$
By the boundedness of $A^{*}$ and the remark \eqref{norm.signs.2}, we obtain:
\begin{align*} 
\sum_{m,n=1}^{M,N}&\left\vert b_{m}^{*}\otimes y_{n}(u_{k})A^{*}(a_{m}\otimes x_{n})\right\vert\\
&=
\left\vert A^{*}\left(\sum_{m,n=1}^{M,N}\varepsilon_{m,n}b_{m}^{*}\otimes y_{n}(u_{k})a_{m}\otimes x_{n}\right)\right\vert\\
&\leq
\Vert A^{*}\Vert_{\left(E\otimes_{\alpha_{\mathcal{F},\mathcal{G}}^{Bess}}F\right )^{*}}
\alpha_{\mathcal{F},\mathcal{G}}^{Bess}\left(\sum_{m,n=1}^{M,N}\varepsilon_{m,n}b_{m}^{*}\otimes y_{n}(u_{k})a_{m}\otimes x_{n}\right)\\
&=
\Vert A^{*}\Vert_{\left(E\otimes_{\alpha_{\mathcal{F},\mathcal{G}}^{Bess}}F\right )^{*}}
\alpha_{\mathcal{F},\mathcal{G}}^{Bess}\left(\sum_{m,n=1}^{M,N}b_{m}^{*}\otimes y_{n}(u_{k})a_{m}\otimes x_{n}\right)\\
&\leq
\Vert A^{*}\Vert_{\left(E\otimes_{\alpha_{\mathcal{F},\mathcal{G}}^{Bess}}F\right )^{*}}
\lim_{M,N\rightarrow \infty}
\alpha_{\mathcal{F},\mathcal{G}}^{Bess}\left(\sum_{m,n=1}^{M,N}b_{m}^{*}\otimes y_{n}(u_{k})a_{m}\otimes x_{n}\right)\\
&=
\Vert A^{*}\Vert_{\left(E\otimes_{\alpha_{\mathcal{F},\mathcal{G}}^{Bess}}F\right )^{*}}
\alpha_{\mathcal{F},\mathcal{G}}^{Bess}\left(u_{k}\right)
\end{align*}
Letting $k$ and after $M,N$ tend to infinity, we obtain:
\begin{align*}
\underset{m,n=1}{\overset{\infty}{\sum }} 
 \left\vert b_{m}^{*}\otimes y_{n}^{*}(u)A^{*}(a_{m}\otimes x_{n})\right \vert  
 &\leq 
 \Vert A^{*}\Vert_{\left(E\otimes_{\alpha_{\mathcal{F},\mathcal{G}}^{Bess}}F\right )^{*}}\alpha_{\mathcal{F},\mathcal{G}}^{Bess}\left(u\right).
\end{align*}
That is, $\mathcal{F}\otimes\mathcal{G}$ is a BSB for 
$E\otimes_{\alpha_{\mathcal{F},\mathcal{G}}^{Bess}}F$.
By theorem \eqref{characterization.USB.Singer},  $\mathcal{F}\otimes\mathcal{G}$ is a USB for $E\otimes_{\alpha_{\mathcal{F},\mathcal{G}}^{Bess}}F$. Thus, the proof is complete.\qed
\begin{proposition}
Suppose that \( E \) and \( F \) are real Banach spaces with unconditional Schauder bases \( \mathcal{F} = ((a_m, b_m^*)) \), and \( \mathcal{G} = ((x_n, y_n^*)) \), respectively. Then the spaces $\left(E,\Vert\cdot\Vert^{Bess}_{E,\mathcal{F}}\right)$, $\left(F,\Vert\cdot\Vert^{Bess}_{F,\mathcal{G}}\right)$ and $\left(E\otimes_{\alpha_{\mathcal{F},\mathcal{G}}^{Bess}}F,\alpha_{\mathcal{F},\mathcal{G}}^{Bess}\right)$ are Banach lattices.
\end{proposition}
\textbf{Proof.}
\begin{enumerate}
\item
Definition of the Partial Order:\\
Since \( \mathcal{F} \otimes \mathcal{G}\) forms a USB for \( E \otimes_{\alpha_{\mathcal{F},\mathcal{G}}^{\text{Bess}}} F \), every \( u \in E \otimes_{\alpha_{\mathcal{F},\mathcal{G}}^{\text{Bess}}} F \) admits a unique expansion:
$$u = \sum_{m,n=1}^{\infty}\lambda_{m,n}^{u} a_m \otimes x_n. $$ 
Define \( u \preceq v\) if and only if \( \lambda_{m,n}^{u}\leq \lambda_{m,n}^{v} \) for all \( m, n \in\mathbb{N}^{*}\).
\item
Verification of Banach Lattice Axioms:
\begin{enumerate}
\item
Translation invariance:\\
 If \( u \preceq  v \), then \( \lambda_{m,n}^u \leq \lambda_{m,n}^v \) for all \( m, n \). For any \( w \), the coefficients of \( u + w \) and \( v + w \) are \( \lambda_{m,n}^u + \lambda_{m,n}^w \) and \( \lambda_{m,n}^v + \lambda_{m,n}^w \), respectively. Thus,
$$\lambda_{m,n}^u + \lambda_{m,n}^w \leq \lambda_{m,n}^v + \lambda_{m,n}^w,$$
which implies:
$$u + w \preceq  v + w.$$
\item
Positive homogeneity:\\
If \(0 \preceq u\) and \(\lambda\geq 0 \), then \( \lambda_{m,n}^{u} \geq 0 \) for all \( m, n \), so \( \lambda \lambda_{m,n}^{u} \geq 0 \), implying 
$$0\preceq \lambda u .$$
\item
 Existence of least upper bound and greatest lower bound:\\
  For any \( u \) and \( v \), define \( w = u \lor v \) by
$$ \lambda_{m,n}^w = \max(\lambda_{m,n}^u, \lambda_{m,n}^v),\;(m,n\in\mathbb{N}^{*}).$$
 Due to the unconditionality of the basis, the series for \( w \) converges in norm, and \( w \) is the least upper bound of \( u \) and \( v \). Similarly, define \( z = u \land v \) by
$$ \lambda_{m,n}^z = \min(\lambda_{m,n}^u, \lambda_{m,n}^v) ,\;(m,n\in\mathbb{N}^{*}).$$
\item
Norm compatibility:\\
The absolute value is defined as
$$ |u| = u \lor (-u),$$
 which has coefficients \( |\lambda_{m,n}^{u}| \). If \( |u| \preceq  |v| \), then \( |\lambda_{m,n}^u| \leq|\lambda_{m,n}^v| \) for all \( m, n \). From proposition \eqref{norm.for.BSB},
   \[
\alpha_{\mathcal{F},\mathcal{G}}^{\text{Bess}}(u)
=
\underset{N\rightarrow\infty}{\lim}
\underset{M\rightarrow\infty}{\lim}
\sup_{\substack{\|f^*\|_{E^*} \leq 1 \\ \|g^*\|_{F^*} \leq 1}} \sum_{m,n=1}^{M,N} |\lambda_{m,n}^{u}| |f^*(a_m)| |g^*(x_n)|,
   \]
it follows that \( \alpha_{\mathcal{F},\mathcal{G}}^{\text{Bess}}(u) \leq \alpha_{\mathcal{F},\mathcal{G}}^{\text{Bess}}(v) \).
Thus, \(\left(E \otimes_{\alpha_{\mathcal{F},\mathcal{G}}^{\text{Bess}}} F,\alpha_{\mathcal{F},\mathcal{G}}^{\text{Bess}}\right) \) is a Banach lattice.\\
The proofs for  $(E,\Vert\cdot\Vert^{Bess}_{E,\mathcal{F}})$ and $(F,\Vert\cdot\Vert^{Bess}_{F,\mathcal{G}})$ are similar.\qed
\end{enumerate}
\end{enumerate}
\section*{Declarations}
\begin{enumerate}
\item 
\textbf{Data Availability Statement.}\\
Not applicable.
\item
\textbf{Competing Interests.}\\
I have nothing to declare.
\item
 \textbf{Funding.}\\
This research did not receive any specific grant from funding agencies in the public, commercial, or not-for-profit sectors.
\item
\textbf{Authors contributions.}\\
The authors equally conceived of the study, participated in its design and coordination, drafted the manuscript, participated in the sequence alignment, 
and read and approved the final manuscript.
\item
\textbf{Ethical approval.}\\
This article does not contain any studies with animals performed by any of the authors.
\end{enumerate}

\bigskip
\bigskip  
\end{document}